\documentclass{article}
\usepackage{graphicx}
\usepackage{amssymb}  
\usepackage{amsthm}
\usepackage{color}
\usepackage{subfloat}
\usepackage{subcaption}
\usepackage{todonotes}
\DeclareMathAlphabet{\mathmybb}{U}{bbold}{m}{n}
\newcommand{\1}{\mathmybb{1}}
\usepackage[shortlabels]{enumitem}

\usepackage[normalem]{ulem}
\usepackage{hyperref}

\usepackage{amsmath}

\newcommand{\ca}{{\bf \vartheta}}
\newcommand{\diageta}{{\bf \Xi}}

\newcommand{\be}{\begin{equation}}
\newcommand{\ee}{\end{equation}}
\newcommand{\bes}{\begin{equation*}}
\newcommand{\ees}{\end{equation*}}
\newcommand{\bee}{\begin{eqnarray*}}
\newcommand{\eee}{\end{eqnarray*}}
\newcommand{\bea}{\begin{eqnarray}}
\newcommand{\eea}{\end{eqnarray}}

\newcommand{\R}{\mathbb{R}}
\newcommand{\B}{\mathcal{B}}

\newtheorem{theorem}{Theorem}[section]

\newtheorem{lemma}[theorem]{Lemma}
\newtheorem{corollary}[theorem]{Corollary}
\newtheorem{remark}[theorem]{Remark}

\newcommand{\comment}[1]{}

\newcommand{\E}{\mathbb E}
\newcounter{algo}[section]

\def\IR{\hbox{\rm I\kern-.2em\hbox{\rm R}}}
\def\IC{\hbox{\rm C\kern-.58em{\raise.53ex\hbox{$\scriptscriptstyle|$}}
    \kern-.55em{\raise.53ex\hbox{$\scriptscriptstyle|$}} }}

\newcommand{\N}{\mathbb{N}}

\newcommand{\xbf}{{x}}
\newcommand{\ybf}{{y}}
\newcommand{\ebf}{\mathbf{e}}
\newcommand{\tbf}{t}

\def \E {{\mathbb E}}

\usepackage{bbm}				
\def \1{{\mathbbm 1}}

\DeclareMathOperator{\diag}{diag}
\DeclareMathOperator{\diam}{diam}
\DeclareMathOperator{\erg}{erg}
\DeclareMathOperator*{\argmin}{arg\,min}
\DeclareMathOperator*{\essinf}{ess\,inf}

\newcommand{\lr}[1]{\left(#1\right)}

\newcommand{\lrq}[1]{\left[#1\right]}

\title{A discrete Consensus-Based Global Optimization Method with Noisy Objective Function}

\author{ Stefania Bellavia\footnotemark[1], Greta Malaspina\footnotemark[1]\footnotemark[1]}

\date{}
\begin {document}

\maketitle
\footnotetext[1]{Dipartimento  di Ingegneria Industriale, Universit\`a degli Studi di Firenze,
Viale G.B. Morgagni 40,  50134 Firenze,  Italia. Members of the INdAM Research Group GNCS. Emails:
stefania.bellavia@unifi.it, greta.malaspina@unifi.it}

\maketitle

\begin{abstract}
\noindent Consensus based optimization is a derivative-free particles-based method for the solution of global optimization problems. Several versions of the method have been proposed in the literature, and different convergence results have been proved. 
However, all existing results assume the objective function to be evaluated exactly at each iteration of the method. In this work, we extend the convergence analysis of a discrete-time CBO method to the case where only a noisy stochastic estimator of the objective function can be computed at a given point. In particular we prove that under suitable assumptions on the oracle's noise, the expected value of the mean squared distance of the particles from the solution can be made arbitrarily small in a finite number of iterations. Numerical experiments showing the impact of noise are also given.
\end{abstract}

\section{Introduction}
We consider a global optimization problem of the form 
\be
x^* = \argmin_{\xbf\in\R^d} f(\xbf)
\ee
where $f:\R^d\longrightarrow\R$ is a smooth, generally non-convex, function which is not available for exact evaluation. That is, we assume that given a point $\xbf\in\R^d$ only a stochastic estimator $\hat f(\xbf,\omega(\xbf)) \simeq f(\xbf)$ con be computed, where 
$\omega(x)$ is a random variable.

Global optimization problems are 
widespread in the mathematical modeling of real problems and a large variety of computational approaches have been proposed for their solution. For an  introduction on the topic, see e.g.  \cite{locatelli2013,horst2000,GObook1,GObook2}.

Consensus Based Optimization (CBO) methods are derivative-free, particle-based optimization methods that employ ideas from drift-diffusion dynamics and consensus formation in multi-agent systems to solve global optimization problems. This class of methods was first proposed in \cite{CBO1} and several papers in the literature study its convergence properties. 
The method considered in \cite{CBO2,Fornasier1,CBO1} is as follows. Consider a set of $N$ particles, where for $i=1,\dots,N$ the vector $\xbf^i_t\in\R^d$ denotes the position of the $i$-th particle at time $t$. Given the initial distribution $\rho_0$, the particles evolve according to the following stochastic differential equation:
\be\label{CBO_continuous}
\begin{cases}
\xbf^i_0\sim\rho_0\\
d\xbf^i_{t} = -\gamma\lr{\xbf^i_t-\xbf^\alpha_t}dt+\sigma\|\xbf^i_t-\xbf^\alpha_t\|_2dB^i_t.
\end{cases}
\ee
where $\gamma \in (0,1]$, $\sigma>0$ are the drift and the diffusion parameter respectively, $\{\{B^i_t\}_{t\geq0}| i=1,\dots,N\}$ are independent Brownian motions in $\R^d$ and for a given value of $\alpha>0$ the consensus point $\xbf^\alpha_t$ is defined as 
{
\bes
\xbf^\alpha_t = \frac{\sum_{i=1}^N \xbf^i_te^{-\alpha f(\xbf^i_t)}}{\sum_{i=1}^N e^{-\alpha f(\xbf^i_t)}}.
\ees}
In \cite{CBO2} the authors perform the theoretical analysis of the method described in \eqref{CBO_continuous} by studying the evolution of its mean-field limit. 
First of all, they prove that, under suitable regularity assumptions on the objective function, the particles achieve consensus. More precisely, they prove that, asymptotically, the distribution of the mean-field limit tends to a distribution centered at one point. Finally, they provide assumptions on the initial distribution of the particles, that ensure that the consensus point is in a neighborhood of the global minimizer of the objective function, where the size of the neighborhood depends on the parameter $\alpha.$
In \cite{Fornasier1} the authors prove that the probability distribution of the mean-field limit of the particles converges to a distribution centered at the solution of the problem, in the Wasserstein distance, with minimal assumptions on the initial distribution of the particles. Moreover, they provide a probabilistic convergence result in the finite-particles regime. Specifically, they prove that the squared norm of the difference between the average particle and the solution becomes smaller than any given threshold with a certain probability. 

A modification of the CBO method given in \cite{CBO1} 
has been proposed in \cite{CBObatches} where the isotropic geometric Brownian motion is replaced by a component-wise one, thus removing the dimensionality dependence of the drift rate.
The stochastic differential equation that describes the evolution of the particle is then given by 
\be\label{CBO_continuous2}
\begin{cases}
\xbf^i_0\sim\rho_0\\
d\xbf^i_{t} = -\gamma\lr{\xbf^i_t-\xbf^\alpha_t}dt+\sigma
\sum_{s=1}^d(x^i_{t,s}-x^\alpha_{t,s})dB_{t,s}^i\ebf_s 
\end{cases}
\ee
where $\ebf_s $ denotes the $s$-th vector of the canonical basis of $\R^d$ 
and for any vector $\xbf\in\R^d$ we denote with $x_s$ the $s$-th component of $\xbf.$
For the finite-sum minimization problem the authors also propose an additional modification that uses a mini-batch strategy on both the particles and the objective function at each iteration. However, the convergence analysis is only carried out in the case where full sample is used to evaluate the objective function.

The following discrete-time version of the method from \cite{CBObatches} with homogeneous noises has been introduced in 
\cite{discreteCBO}: 
\be\label{CBO_ft}
\begin{cases}
\xbf^i_0\sim\rho_0\\
\xbf^i_{k+1} = \xbf^i_k-\gamma\lr{\xbf^i_k-\xbf_k^\alpha}-\sum_{s=1}^d\lr{x^i_{k,s}- x^\alpha_{k,s}}\eta_{k,s}\ebf_s,
\end{cases}
\ee
where $\{\eta_{k,s}\}$ are i.i.d. random variables with mean zero and variance $\sigma^2$, {and the consensus point $\xbf^\alpha_k$ is defined as 
\be\label{xkalpha}
\xbf^\alpha_k = \frac{\sum_{i=1}^N \xbf^i_ke^{-\alpha f(\xbf^i_k)}}{\sum_{i=1}^N e^{-\alpha f(\xbf^i_k)}}.
\ee}
Finally, \cite{hetCBO} analyses a modification of \eqref{CBO_ft} that, analogously to \cite{CBObatches}, employs heterogeneous noise and random batch interaction for the computation of the consensus point $\xbf^\alpha_k.$ The authors prove that under suitable choices of the parameters of the method the consensus is formed with probability one.

Several extensions of the CBO method have been proposed in literature for the solution of different classes of optimization problems that ranges from 
the optimization of non-convex functions on compact hypersurfaces \cite{CBOhyper} and 
application of the method to machine learning problems \cite{CBOsphere} to 
the solution of constrained optimization problems \cite{CBOconstr} and 
multi-objective optimization  \cite{CBO_MOO}.

The theoretical analysis of all the methods mentioned so far assumes that the objective function is evaluated exactly at every iteration. However, in many practical applications the exact value of the objective function is not available.
The noise might for instance stem from propagation of round-off
errors or from  numerical simulations that provides only an approximation of the objective function \cite{parsopoulos2002recent,sergeyev2020safe,villemonteix2009informational}. Furthermore, in many
problems the accurate values of the function to be minimized are computationally expensive to obtain and have to be approximated by stochastic estimators, as 
in the case of finite-sum minimization problems, see e.g   the empirical-risk minimization in machine learning models training \cite{candelieri2019global,CBObatches}.
In this paper, we provide a theoretical analysis of CBO method in the case of noisy objective functions. On the one hand, we prove that the theory presented in \cite{discreteCBO} can be extended to the case of noisy objective function and {node-dependent diffusion terms $\eta^i_{k,s}$}, given reasonable assumptions on the stochastic estimator. On the other hand, under regularity assumptions that are analogous to those in \cite{Fornasier1} and suitable assumptions on the stochastic estimator and on the particles position, we provide a  worst-case iteration complexity result in expectation, showing that, given $\epsilon>0$ the expected value of the average distance of the particles to the global minimizer is smaller than $\epsilon$ after at most $O(\log(\epsilon^{-1}))$ iterations.
To the best of our knowledge the second result, namely, the complexity and the consequent convergence result to zero of the average distance from the solution, is not available in literature for discrete CBO, even when the objective function is evaluated exactly. Then, our results also complement those in \cite{discreteCBO,hetCBO}.

Some numerical experiments are also presented showing the impact of noise.
We consider both the case of additive Gaussian noise and of noise arising by subsampling estimates of the objective function in the context  of binary classification problems.
The results show that subsampling can significantly reduce the overall computational cost of the method while retaining a good level of accuracy.

This paper is organized as follows. In Section \ref{sec:method} we describe the proposed approach and in Section \ref{sec:analysis} we present its convergence analysis. 
Section \ref{sec:numres} shows  some computational evidence of the reliability and robustness to noise of CBO methods.
{In Section \ref{sec:conclusions} we provide some conclusions.}

{\it Notations.} Throughout the paper, $\R$ and $\N$ denote the set of real  and natural numbers, respectively and $\|\cdot\|$ denotes the 2-norm.  Given a random variable $X:\Omega\rightarrow \R$, we denote with $\E[X]$ the expected value of $X$.
We will further denote with $\otimes$ the Kronecker product and given an integer $m$ with $\1_m$ the vector of all ones of dimension $m$  and with $I_m$ the $m\times m$ identity matrix.

\section{Consensus based optimization with noisy functions }\label{sec:method}
Let us assume that we have $N\in\N$ particles, and for $i=1,\dots,N$, $k\in\N$ let us denote with $\xbf^i_k$ the position of particle $i$ at iteration $k$. 
We further assume that we have at disposal only noisy evaluations of the objective function at each particle $i$. 
More precisely, 
the function's estimates are produced by a stochastic oracle such that, given a point $\xbf$  computes a function estimate 
$\hat f(\xbf, \omega(\xbf))$
where $\omega(\xbf)$ is the randomness of 
the oracle and depends on the input value $\xbf$. 
For sake of brevity, given the generic  particle $\xbf^i_k$
we will denote the oracle's estimate at $\xbf^i_k$  as 
\be\label{noisy_f}
\hat f(\xbf^i_k, \omega^i_k).
\ee
We consider the discrete analogue of \eqref{CBO_continuous2} and,
for every $i$, the sequence $\{\xbf^i_k\}$ is then generated as follows

\be\label{CBOhat}
\begin{cases}
\xbf^i_0\sim\rho_0\\
\xbf^i_{k+1} = \xbf^i_k-\gamma\lr{\xbf^i_k-\hat\xbf^\alpha_k}-\sum_{s=1}^d\lr{x^i_{k,s}-\hat x^\alpha_{k,s}}\eta_{k,s}^i\ebf_s
\end{cases}
\ee
where $\gamma>0$ is the drift parameter, $x^i_{k,s}$ and $\hat x^\alpha_{k,s}$ denote the component $s$ of $x^i_{k}$ and $\hat x^\alpha_{k}$ respectively, 
$\{\eta_{k,s}^i\}_{s,k,i}$ are i.i.d. Gaussian random variables with 
$$
\E[\eta_{k,s}^i] = 0 \quad \E[(\eta_{k,s}^i)^2] = \xi^2,$$
for some $\xi>0$, and for a given value of $\alpha>0$ the consensus point $\hat\xbf^\alpha_k$, that depends on the noisy $f$-evaluations, is defined as
\be \label{noisy_cp}
\hat\xbf^\alpha_k = \frac{\sum_{i=1}^N \xbf^i_ke^{-\alpha\hat f(\xbf^i_k, \omega^i_k)}}{\sum_{i=1}^N e^{-\alpha\hat f(\xbf^i_k, \omega^i_k)}}.
\ee

Notice that by Lemma \ref{lemma:gbounds} in the Appendix for any $s=1\ldots, d$ and any $k$ it holds
\begin{equation}\label{expmax}
    \begin{aligned}
    &\E\lrq{\max_{i=1:N}|\eta^i_{k,s}|}\leq 2\xi(\log(\sqrt{2}N))^{1/2},\enspace \enspace\E\lrq{\max_{i=1:N}|\eta^i_{k,s}|^2}\leq 4\xi^2\log(\sqrt{2}N).
    \end{aligned}
\end{equation}
We also introduce the constant 
\be \label{a}
\ca = 1-\gamma+8\xi(\log(\sqrt{2}N))^{1/2}
\ee
that depends on the number of particles $N$ and the parameters $\gamma,\xi$.
Let $\diageta_k^i\in \R^{d\times d} $ the diagonal matrix  
such that $\diageta_k^i=\diag(\eta_{k,1}^i,\ldots,\eta_{k,d}^i)$
  and 
$$
(\Delta \xbf)_k^i=(\gamma I_d+\diageta_{k}^i)\lr{\hat x^\alpha_{k}-x^i_{k}}.
$$
Then, it holds
\be \label{CBOhat_compatto}
\xbf^i_{k+1} = \xbf^i_k+(\Delta \xbf )_k^i.
\ee
Following \cite{hetCBO} we will also make use of the following matrix reformulation of the evolution of each component of the particles.
For any $s=1,\dots,d$ and $k\in\N$ let us define $\ybf_k^s = \lr{x^1_{k,s},\dots,x^N_{k,s}}^\top\in\R^N$ and $D^s_k $ the diagonal matrix such that  $D^s_k = \diag(\eta^1_{k,s},\dots,\eta^N_{k,s})\in\R^{N\times N}$.  By \eqref{CBOhat} we have

\begin{equation}\label{Y1}
    \ybf^s_{k+1} = \ybf^s_{k} - \gamma(\ybf^s_{k} - \hat\xbf^\alpha_{k,s}\otimes \1_N ) - D^s_k (\ybf^s_{k} - \hat\xbf^\alpha_{k,s}\otimes \1_N ).
\end{equation}
By definition of $\hat\xbf^\alpha_k$ we further have that for any $s=1,\dots,d$ and $k\in\N$
\be \label{xhat}
\hat\xbf^\alpha_{k,s} = \sum_{i=1}^N \hat v^\alpha_{k,i}\xbf^i_{k,s}\enspace \text{with}\ \hat v^\alpha_{k,i}  = \frac{e^{-\alpha\hat f(\xbf^i_k, \omega^i_k)}}{\sum_{i=1}^N e^{-\alpha\hat f(\xbf^i_k, \omega^i_k)}}
\ee
and therefore
$\hat\xbf^\alpha_{k,s}\otimes \1_N = \hat V^\alpha_k \ybf^s_k,$
where $\hat V^\alpha_k\in\R^{N\times N}$ has all rows equal to $(\hat v^\alpha_{k,1}, \dots, \hat v^\alpha_{k,N} )$.
Equation \eqref{Y1} can then be rewritten as
\begin{equation}
\begin{aligned}\label{Y2}
    \ybf^s_{k+1}& = \ybf^s_{k} - \gamma(I_N-\hat V^\alpha_k)\ybf^s_k - D^s_k (I_N-\hat V^\alpha_k)\ybf^s_k \\
    &= \lr{(1-\gamma)I_N + \gamma \hat V^\alpha_k  - D^s_k (I_N-\hat V^\alpha_k)}\ybf^s_k =: M^s_k\ybf^s_k.
\end{aligned}
\end{equation}

\section{Convergence Analysis}\label{sec:analysis}
In the subsequent analysis we will make use of the following assumptions.\\
\noindent\textbf{Assumption A1}. The objective function $f:\R^d\longrightarrow \R$ is twice-continuously differentiable and strictly positive.

\noindent\textbf{Assumption A2}. There exists $M_H>0$ such that $\|\nabla^2 f(\xbf)\|\leq M_H$ for every $\xbf\in\R^d.$

\noindent\textbf{Assumption A3}. There exists a unique $\xbf^*\in\R^d$ such that 
$$\xbf^* = \argmin_{\xbf \in \R^d}f(x).$$\\
Moreover, $f$ is locally strictly convex at $\xbf^*$. That is, $\nabla^2f(\xbf^*) $ is positive definite.
In the following we define
$f_* = f(\xbf^*).$

\noindent\textbf{Assumption A4}. For every $i=1,\dots,N$ the initial position $\xbf^i_0\in\R^d$ is a random variable with law $\rho_0$, 
where $\rho_0$ is absolutely continuous with respect to the Lebesgue measure. Moreover, the associated density function $\varphi_0$ has compact support, is continuous at $\xbf^*$ and satisfies $\varphi_0(\xbf^*)>0.$

\noindent\textbf{Assumption A5}. The noisy oracle $\hat f$ is a function of $\xbf\in\R^d$ and of a random variable $\omega(\xbf)$ that depends on $\xbf$. Moreover, there exist two constants $t_0,\ t_1\geq 0$ such that 
\be\label{bound_noise}
\E_{\omega(\xbf)}[|f(\xbf)-\hat f(\xbf,\omega(\xbf))|^2\ ]\leq t_0+t_1 f(\xbf)^2\enspace\forall\xbf\in\R^d,
\ee
where $E_{\omega(\xbf)}$ denotes the expected value with respect to the distribution of the random variable 
$\omega(\xbf)$ given $\xbf$.\\

\begin{remark}
Assumption A4 requires the unknown solution $x^*$ to be in the support of the initial distribution $\rho_0$ of the particles. This can be achieved, in practice, taking a distribution $\rho_0$ with support over the whole space $\R^n$. This assumption is necessary for our theoretical analysis in the discrete-time finite-particles regime. In \cite{Fornasier2} the authors prove, in the continuous-time setting, that the convergence of CBO method to the global minimizer can be proved also when the initial distribution has no mass around $x^*.$
\end{remark}

\begin{remark}
Assumption A5 does not require the variance of the oracle to go to zero and  also applies to the setting in which available function value  estimates are occasionally biased. If the oracle is unbiased,\eqref{bound_noise} is equivalent to
\be\label{bound_noise_unbiased}
\E_{\omega(\xbf)}[\hat f(\xbf,\omega(\xbf))^2\ ]\leq t_0+(t_1+1) f(\xbf)^2\enspace\forall\xbf\in\R^d.
\ee 
Then we assume that the variance of the function estimate satisfies a growth condition. This is weaker and more realistic than the bounded variance condition and it is  used in the analysis of  stochastic optimization methods, see e.g \cite{bottou,fang2024fully}. Note that in particular Assumption A5 does not require the expectation of the squared error to vanish as the true value of the function approaches 0.
\end{remark}

We start  by studying the asymptotic behavior of the generated sequences $\{x_k^i\}$, $i=1,\ldots,N$. Namely, we prove that with probability one, all the particles converge to the same point $\xbf_{\infty}$, as $k$ tends to $+\infty.$ The proof of convergence follows the path of \cite{hetCBO}.
We first provide some auxiliary results.

\subsection{Diameter and Ergodicity}

Given any vector $\ybf\in\R^N$ and any matrix $M\in\R^{N\times N}$ we denote with $\diam(\ybf)$ and $\erg(M)$ the diameter of $\ybf$ and the ergodicity coefficient of $M$, respectively. That is, we define
\begin{eqnarray} 
\diam(\ybf) &=& \max_{i,l = 1:N} y_i-y_l = \max_{i=1:N}y_i - \min_{i=1:N}y_i={\max_{i,l=1:N}} |y_i-y_l| \label{diam} \\
\erg(M) &=& \min_{i,l=1:N}\sum_{j=1}^N \min\{M_{ij}, M_{lj}\}. \label{ergo}
\end{eqnarray}
Notice that $\diam(\ybf) = 0$ if and only if $y_i =y_l$, for any $i,l=1,\ldots,N$. We
will now  provide upper bounds for the evolution of  $\diam(\ybf^s_k)$, as well as for the maximum distance between any two particles, and the expected  distance between a particle and  the consensus point $\hat\xbf^\alpha_k$. We will make use of the following   preliminary result
whose proof is provided in the appendix.

\begin{lemma}\label{lemma_ergM}
Let $\ybf^s_{k+1}$  and $M^s_k$ be given  in \eqref{Y2}.
For any $k\in\N$ and $s=1,\dots,d$ we have
\begin{description}
    \item{i)} $\diam(\ybf^s_{k+1})\leq(1-\erg(M_k))\diam(\ybf^s_k)$
      \item{ii)} $\erg(M^s_k)\geq \gamma -4\max_{i=1:N}|\eta^i_{k,s}|.$
\end{description}
\end{lemma}

\begin{lemma}\label{lemma_diam_CBO}Let $\xbf_k^i$, $i=1,\ldots N$ be the $N$ particles generated at iteration $k$ of \eqref{CBOhat} and  $\ybf_k$ and $\ca$ given in \eqref{Y1} and \eqref{a}, respectively.
For every $s=1,\dots,d$ and $k\in\N$
\begin{enumerate}[i)]
\item $\E[\diam(\ybf^s_{k})]\leq \ca^k\E[\diam(\ybf^s_{0})]$
\item $\diam(\ybf^s_{k})\leq e^{-kS_{k,s}}\diam(\ybf^s_{0})$
where $S_{k,s}$ is a random variable such that 
$$\lim_{k\rightarrow +\infty}S_{k,s} \ge 1-\ca\enspace \text{almost surely}$$ 
\item $\displaystyle\max_{i,j=1:N}\|\xbf^i_{k}-\xbf^j_{k}\|_\infty\leq \max_{s=1:d} e^{-kS_{k,s}}\max_{i,j=1:N}\|\xbf^i_{0}-\xbf^j_{0}\|_\infty$
\item $\E\lrq{\|\xbf^i_k - \hat\xbf^\alpha_k\|^2}\leq e^{-2k(1-\ca)}\sum_{s=1}^d\E\lrq{\diam(\ybf^s_0)^2}$.

\end{enumerate}
\end{lemma}
\begin{proof}
Lemma \ref{lemma_ergM} yields
\begin{equation}\label{diam_rec}
\begin{aligned}
\diam\lr{\ybf^s_{k}}&\leq\lr{1-\gamma+4\max_{i=1:N}|\eta^i_{k-1,s}|}\diam\lr{\ybf^s_{k-1}}\\
&\leq \diam\lr{\ybf^s_{0}}\prod_{l=0}^{k-1}\lr{1-\gamma+4\max_{i=1:N}|\eta^i_{l,s}|}.
\end{aligned}
\end{equation}
By \eqref{expmax}, we have
\begin{equation}\label{prodexp}
\begin{aligned}
    \E\lrq{\prod_{l=0}^{k-1}\lr{1-\gamma+4\max_{i=1:N}|\eta^i_{l,s}|}} &= \prod_{l=0}^{k-1}\lr{1-\gamma+4\E\lrq{\max_{i=1:N}|\eta^i_{l,s}|}}\leq\ca^k.
\end{aligned}
\end{equation}
Taking the expectation in \eqref{diam_rec} and using this inequality we  get $i)$.\\
To prove $ii)$ let us first notice that, since $e^t\geq t+1$ for every $t\in\R$, we get 
$1-\gamma+4\max_{i=1:N}|\eta^i_{l,s}|\le e^{-\gamma+4\max_{i=1:N}|\eta^i_{l,s}|}$ and 
from \eqref{diam_rec}
\begin{equation}\label{diam_rec_exp}
\begin{aligned}
\diam\lr{\ybf^s_{k}}&\leq \diam\lr{\ybf^s_{0}}\exp\lr{\sum_{l=0}^{k-1}\lr{-\gamma+4\max_{i=1:N}|\eta^i_{l,s}|}}\\
& = \diam\lr{\ybf^s_{0}}\exp(-kS_{k,s})
\end{aligned}
\end{equation}
with $S_{k,s}$ random variable given by 
$$S_{k,s} = \frac{1}{k}\sum_{l=0}^{k-1}\lr{\gamma - 4\max_{i=1:N}|\eta^i_{l,s}|}.$$
Since $\{\eta^i_{0,s},\dots,\eta^i_{k,s}\}$ are i.i.d., using the law of large numbers and \eqref{expmax} we get 
$$\lim_{k\rightarrow+\infty}S_{k,s} = \E\lrq{\gamma - 4\max_{i=1:N}|\eta^i_{l,s}|} \geq\gamma-8\xi\log(\sqrt{2}N)^{1/2}=1-\ca\enspace \text{almost surely},$$ 
and therefore $ii)$ is proved.
By part $ii)$ and the definition of $\ybf^s$ we have
\begin{eqnarray*}
    \max_{i,j=1:N}\|\xbf^i_{k}-\xbf^j_{k}\|_\infty &=&\max_{i,j=1:N}\max_{s=1:d} |x^i_{k,s}-x^j_{k,s}| = \max_{s=1:d} \diam(\ybf^s_k)\\
    &\leq& \max_{s=1:d}e^{-kS_{k,s}} \diam(\ybf^s_0)
    \le\max_{s=1:d}e^{-kS_{k,s}} \max_{s=1:d}\max_{i,j=1:N} |x^i_{0,s}-x^j_{0,s}|\\
    &=& \max_{s=1:d}e^{-kS_{k,s}}\max_{i,j=1:N}\|\xbf^i_{0}-\xbf^j_{0}\|_\infty.
\end{eqnarray*}
which is $iii).$
We now prove $iv)$. Let us first notice that using 
\eqref{xhat},
the fact that $\sum_{j=1}^{N}\hat v^\alpha_{k,j} = 1$ and \eqref{diam_rec} we have
    \begin{eqnarray}
    |x^i_{k,s}-\hat x^\alpha_{k,s}| &=&\left|x^i_{k,s} - \sum_{j=1}^N \hat v^\alpha_{k,j}x^j_{k,s}\right| 
    \leq \sum_{j=1}^N \hat v^\alpha_{k,j}|x^i_{k,s} - x^j_{k,s}|\leq \max_{j=1:N}|x^i_{k,s} - x^j_{k,s}|\nonumber \\ 
   & =&\diam(\ybf^s_k) \leq \diam\lr{\ybf^s_{0}}\prod_{l=0}^{k-1}\lr{1-\gamma+4\max_{i=1:N}|\eta^i_{l,s}|}.
   \label{th2_1}
    \end{eqnarray}
    Therefore, \eqref{expmax} and \eqref{a} yield
\be \label{bound1}
    \begin{aligned}
    \E\lrq{(x^i_{k,s}-\hat x^\alpha_{k,s})^2}& \leq \E\lrq{\diam\lr{\ybf^s_{0}}^2}\prod_{l=0}^{k-1}\E\lrq{\lr{1-\gamma+4\max_{i=1:N}|\eta^i_{l,s}|}^2} \\&
    \leq 
\E\lrq{\diam\lr{\ybf^s_{0}}^2}\ca^{2k}.
    \end{aligned}
    \ee
Using the fact that $e^{\ca-1}\geq \ca$  we then have the thesis:
\begin{equation*}
    \begin{aligned}
        \E\lrq{\|\xbf^i_k-\hat\xbf^\alpha_k\|^2}& = \sum_{s=1}^d\E\lrq{(x^i_{k,s}-\hat x^\alpha_{k,s})^2}\leq \ca^{2k}\sum_{s=1}^d \E\lrq{\diam\lr{\ybf^s_{0}}^2}\\
&\leq e^{-2k(1-\ca)}\sum_{s=1}^d\E\lrq{\diam(\ybf^s_0)^2}.
    \end{aligned}
\end{equation*}
\qed\end{proof}

\subsection{Almost sure consensus }
In this subsection, assuming $\ca<1$ and
exploiting Lemma \ref{lemma_diam_CBO}, we can show that the particles achieve asymptotic consensus with probability 1. 
\begin{theorem}\label{thm_cons1}
Let $\{(\xbf^i_k)_{i=1:N}\}_{k}$ be the sequence generated by \eqref{CBOhat},$\ca$ given in \eqref{a}, and $\ybf_k^s$,  $s=1,\ldots,d$ given in \eqref{Y1}.
If $\gamma\in(0,1)$, $\xi>0$ and $N\in\N$ are such that $\ca<1$, then
\begin{enumerate}[i)]
\item $\lim_{k\rightarrow+\infty}\E[\diam(\ybf^s_{k})]= 0$;
\item 
$\lim_{k\rightarrow+\infty}\diam(\ybf^s_{k}) = 0$ almost surely;
\item $\lim_{k\rightarrow+\infty}\max_{i,j=1:N}\|\xbf^i_{k+1}-\xbf^j_{k+1}\|_\infty=0$ {almost surely}
\end{enumerate}
\end{theorem}
\begin{proof}
    By Lemma \ref{lemma_diam_CBO} we have 
$ \E[\diam(\ybf^s_{k})] = \mathcal{O}\lr{\ca^k}$, $ \diam(\ybf^s_{k}) = \mathcal{O}\lr{e^{-k(1-\ca)}}$ and
             $\max_{i,j=1:N}\|\xbf^i_{k}-\xbf^j_{k}\|_\infty  = \mathcal{O}\lr{e^{-k(1-\ca)}}$ {almost surely.}
    Then, the thesis follows immediately since $\ca<1$.
\qed\end{proof}

Under the same assumptions of Theorem \ref{thm_cons1}, the following theorem ensures almost sure convergence of the sequence generated by each particle to the same point. The proof, which is analogous to that of Theorem 3.1 in \cite{discreteCBO}, is provided in the Appendix for completeness.
\begin{theorem}\cite{discreteCBO}\label{convergencethm}
    Let $\{(\xbf^i_k)_{i=1:N}\}_{k}$ be the sequence generated by \eqref{CBOhat}. 
    If $\gamma\in(0,1)$, $\xi>0$ and $N\in\N$ are such that $\ca<1$, then
    there exists $\xbf_{\infty}\in\R^d$ such that 
    $$\lim_{k\rightarrow+\infty }\xbf^i_k = \xbf_{\infty}\enspace \text{almost surely},\enspace \text{for every}\ i=1,\dots,N.$$
\end{theorem}

\begin{corollary}\label{boundedfgf}
Let the assumptions in Theorem \ref{convergencethm} hold.
    Then,
    \begin{enumerate}[i)]
    \item there exists a constant $M_f>0$ such that, for every $i=1,\dots,N$ and $k\in\N$, $f(\xbf^i_k)\leq M_f$ almost surely;
    \item $\lim_{k\rightarrow+\infty}\hat\xbf^\alpha_k = \xbf_\infty$ almost surely;
      \item there exists a constant $M_g>0$ such that, for every $k\in\N$, $\|\nabla f(\hat\xbf^\alpha_k)\|\leq M_g$ almost surely.
    \end{enumerate}
    \end{corollary}
\begin{proof}
Part $i)$ follows from the continuity of $f$ by noticing that Theorem \ref{convergencethm} implies in particular that, with probability 1, the sequence $\{\xbf^i_k\}_{k=0}^\infty$ is bounded.\\
Part $ii)$ follows directly from the almost sure convergence of $\{\xbf^i_k\}$ for every $i$.\\
Part $iii)$ follows from the fact that $\nabla f$ is continuous and, by $ii)$, the sequence $\{\hat\xbf^\alpha_k\}_{k=0}^\infty$ is bounded. with probability 1.
\qed\end{proof}
\subsection{Handling the noise}
Let us denote with ${\cal E}_k$ the vector whose entries are given by  the absolute error in the evaluation of the objective function at the point $x_k^i$, i.e.
\be \label{errore}
{\cal E}_{k,i}=|f(\xbf^i_k)-\hat f(\xbf^i_k, \omega^i_k)|,\quad i=1\ldots,N.
\ee
and with
    $\sigma_k$  the $\sigma$-algebra generated by $\{\xbf^i_{l,s}\}_{l=0,\ldots,k}$ and 
    $\{\eta_{l,s}\}_{l=0,\ldots,k}$,
    for $i=1,\ldots,N$, $s=1,\ldots,d$.
       The following Lemma bounds the difference between $\xbf^\alpha_k$ and $\hat\xbf^\alpha_k$ in terms of the error in the function evaluation, where
$\xbf^\alpha_k$ is given in \eqref{xkalpha} and is the analogous of $\hat\xbf^\alpha_k$ computed with the exact objective function. 
       
\begin{lemma}\label{lemmadiffxalpha}
Let $\hat  \xbf^\alpha_k$ and  $\xbf^\alpha_k$ be given in \eqref{xhat} and \eqref{xkalpha}, respectively. Then,
    $$|x^\alpha_{k,s}-\hat x^\alpha_{k,s}|\leq \alpha  \|{\cal E}_k\|_{\infty} \diam(\ybf^s_k)\quad \quad \forall s=1,\dots,d,$$
    where ${\cal E}_k$ is given in \eqref{errore}.
\end{lemma}
\begin{proof}
    Given $x^i_{k,s}\in\R$, $i=1,\ldots,N$, $\alpha>0$ we define the function $\Psi_{k,s}:\R^N\longrightarrow\R$  
    $$\Psi_{k,s}(t):=\frac{\sum_{j=1}^Nx^j_{k,s}e^{-\alpha t_j}}{\sum_{j=1}^N e^{-\alpha t_j}},$$
    where $t_j$ is the $j$-th component of the vector $t$.
    $\Psi_{k,s}$ is continuously differentiable and 
$$\frac{\partial\Psi_{k,s}(t)}{\partial t_i} = \alpha e^{-\alpha t_i}\frac{\sum_{j=1}^N e^{-\alpha t_j}(x^j_{k,s}-x^i_{k,s})}{\lr{\sum_{j=1}^N e^{-\alpha t_j}}^2}.$$
    Given $\tbf_{k} = (f(\xbf^1_k),\dots,f(\xbf^N_k))$ and $\hat \tbf_{k} = (\hat f(\xbf^1_k, \omega^1_k),\dots,\hat f(\xbf^N_k,, \omega^N_k))$ we have
    $x^\alpha_{k,s} = \Psi_{k,s}(\tbf_{k})$,  $\hat x^\alpha_{k,s} = \Psi_{k,s}(\hat \tbf_{k})$
    and $|\tbf_{k}-\hat \tbf_{k}|={\cal E}_k$.
    By the mean value theorem for $\tau = \theta\tbf_k+(1-\theta)\hat\tbf_k$ for some $\theta\in(0,1)$, it holds
    $$
    \begin{aligned}
    &|x^\alpha_{k,s} - \hat x^\alpha_{k,s} | = |\Psi_{k,s}(\tbf_{k}) - \Psi_{k,s}(\hat \tbf_{k})| = \nabla\Psi_{k,s}(\tau)^\top(\tbf_k-\hat \tbf_k)  \\  
    &\leq \max_{i=1:N} {\cal E}_{k,i}\sum_{i=1}^N \left|\frac{\partial\Psi_{k,s}(\tau)}{\partial t_i} \right|    = \|{\cal E}_{k}\|_{\infty} \sum_{i=1}^N \alpha e^{-\alpha \tau_i}\frac{\sum_{j=1}^N e^{-\alpha \tau_j}|x^j_{k,s}-x^i_{k,s}|}{\lr{\sum_{j=1}^N e^{-\alpha \tau_j}}^2}  \\  
    &\leq\alpha  \|{\cal E}_{k}\|_{\infty}\diam(\ybf^s_k),
    \end{aligned}
    $$
    which concludes the proof.
\qed\end{proof}
Under Assumption A5 we can now bound the expected values of  $\|{\cal E}_{k}\|_{\infty}$
and  $\max_{i=1:N} {\cal E}_{k,i}^2$.

\begin{lemma}\label{exp_max_err}
Let  $x_k^i$, $i=1,\ldots N$ be the $N$ particles generated at iteration $k$ of \eqref{CBOhat} and 
${\cal E}_{k}$ given in \eqref{errore}.
If Assumptions A1 and A5 hold, then there exists a constant $M_v$ such that $\E\lrq{\|{\cal E}_{k}\|_{\infty}^2}\leq M_v^2$ for every $k\in\N,$ with $M_v= N^{1/2}\lr{t_0+t_1 M_f^2}^{1/2}.$ \\
In addition, if the oracle $\hat f$ satisfies: 
\begin{equation}\label{Gaussian_noise}
\widehat f(\xbf,\omega(\xbf)) = f(\xbf)+\omega_0(\xbf)+ \omega_1(\xbf)f(\xbf),
\end{equation}
with  $\omega(\xbf)\in\R^2$ random vector with Gaussian components $\omega_j(\xbf)\sim\mathcal{N}(0,t_j)$ for $j=0,1$, then $M_v=  2 (t_0 +t_1 M_f^2)^{1/2} \log{(\sqrt{2}N)}^{1/2}.$

\end{lemma}
\begin{proof}
    For any vector $z\in\R^N$ it holds
 $ \|z\|_\infty^2\leq\|z\|_2^2 = \sum_{i=1}^Nz_i^2. $
        Then,
    \begin{equation}\label{exp_max_err1} 
         \E\lrq{\|{\cal E}_{k}\|_{\infty}^2}\leq \E\lrq{{\sum_{i=1}^N {\cal E}_{k,i}^2}}
            \leq \sum_{i=1}^N\E\lrq{{\cal E}_{k,i}^2}.
      \end{equation}
    Let us now consider the $i$-th term of the sum on the right hand side. Using the law of total probability, Assumption A5  and Corollary \ref{boundedfgf}, we have 
    \begin{equation*}
   \E\lrq{{\cal E}_{k,i}^2}= \E\lrq{\E\lrq{{\cal E}_{k,i}^2}|\sigma_k}
    = \E\lrq{t_0+t_1f(\xbf^i_k)^2} = t_0+t_1\E[f(\xbf^i_k)^2]\leq  t_0+t_1 M_f^2.
    \end{equation*}
    Using this inequality in \eqref{exp_max_err1} we get 
    \begin{equation*}
   \E\lrq{\|{\cal E}_{k}\|_{\infty}^2}\leq 
     \sum_{i=1}^N \lr{t_0+t_1 M_f^2} = N\lr{t_0+t_1 M_f^2 }.
    \end{equation*}
    which yields the thesis.  
    Let us now consider the case of additive Gaussian noise. Using \eqref{Gaussian_noise}, the tower rule, and Lemma \ref{lemma:gbounds} in the Appendix we have
    \begin{equation}\begin{aligned}
        \E\lrq{\|\mathcal{E}_k\|^2_{\infty}} &= \E\lrq{\max_{i=1:N}|\omega_0(x^i_k)+\omega_1(x^i_k)f(x^i_k)|^2} \\& =\E\lrq{\E\lrq{\max_{i=1:N}|\omega_0(x^i_k)+\omega_1(x^i_k)f(x^i_k)|^2\ |\ \sigma_k}}\\ &\leq \E\lrq{4\log\lr{\sqrt{2}N}\max_{i=1:N}(t_0+t_1f(x^i_k))},
    \end{aligned}\end{equation}
    where we used the fact that, conditioned to $\sigma_k$, each $\omega_0(x^i_k)+\omega_1(x^i_k)f(x^i_k)$ is a Gaussian random variable with mean 0 and variance $t_0+t_1f(x^i_k)^2.$ The thesis now follows using part $i)$ of Corollary \ref{boundedfgf} and the law of total expectation.

\qed\end{proof}

\subsection{Almost sure convergence}
The results in this subsection prove that, under suitable assumptions on the parameters of the method and the initial distribution of the particles, the optimal gap at the consensus point, namely $|\essinf f(\xbf_{\infty}) - f_*|$, is bounded by a quantity that depends on the dimension of the problem $d$ and the parameter $\alpha$.  To this end we need to provide a lower bound to  the expected value of $\frac{1}{N}\sum_{i=1}^N(e^{-\alpha f(\xbf^i_{k+1})}- e^{-\alpha f(\xbf^i_{k})})$ through the following intermediate results.

\begin{lemma}\label{lemma14}
Let  $x_k^i$, $i=1,\ldots N$ be the $N$ particles generated at iteration $k$ of \eqref{CBOhat} and assume that Assumption A1 holds. Then, 
for every $k\in\N$ there exist $\theta^1_k,\dots,\theta^N_k\in(0,1)$ such that
    \be\label{meandiffexp}
    \frac{1}{N}\sum_{i=1}^Ne^{-\alpha f(\xbf^i_{k+1})}-\frac{1}{N}\sum_{i=1}^N e^{-\alpha f(\xbf^i_{k})}\geq -\alpha A_k-\alpha B_k,
    \ee
    with
    \bes
    \begin{aligned}
    &A_k = \frac{1}{N}\sum_{i=1}^Ne^{-\alpha f(\xbf^i_{k})}\lr{\nabla f(\theta^i_k\xbf^i_{k+1}+(1-\theta^i_k)\xbf^i_k)-\nabla f(\hat\xbf^\alpha_k)}^\top(\xbf^i_{k+1}-\xbf^i_k)\\
    &B_k = \frac{1}{N}\sum_{i=1}^Ne^{-\alpha f(\xbf^i_{k})}\nabla f(\hat\xbf^\alpha_k)^\top(\xbf^i_{k+1}-\xbf^i_k).
    \end{aligned}
    \ees
\end{lemma}
\begin{proof}
Using the fact that $e^y-1\geq y$ for any $y\in\R$ and the mean value theorem, we have that for every $i=1,\dots,N$ 
    \be\label{diffexp}
    \begin{aligned}
    e^{-\alpha f(\xbf^i_{k+1})}-e^{-\alpha f(\xbf^i_{k})}&= e^{-\alpha f(\xbf^i_{k})}\lr{e^{-\alpha( f(\xbf^i_{k+1})-f(\xbf^i_{k}))}-1}   \\  
    & \geq -\alpha e^{-\alpha f(\xbf^i_{k})}\lr{f(\xbf^i_{k+1})-f(\xbf^i_{k})}   \\  
    &=-\alpha e^{-\alpha f(\xbf^i_{k})}\nabla f(\theta^i_k\xbf^i_{k+1}+(1-\theta^i_k)\xbf^i_k)^\top(\xbf^i_{k+1}-\xbf^i_k)   \\  
    &=-\alpha e^{-\alpha f(\xbf^i_{k})}\lr{\nabla f(\theta^i_k\xbf^i_{k+1}+(1-\theta^i_k)\xbf^i_k)-\nabla f(\hat\xbf^\alpha_k)}^\top(\xbf^i_{k+1}-\xbf^i_k)   \\  
    &\quad -\alpha e^{-\alpha f(\xbf^i_{k})}\nabla f(\hat\xbf^\alpha_k)^\top(\xbf^i_{k+1}-\xbf^i_k)
    \end{aligned}
    \ee
    for some $\theta^i_k\in(0,1).$ Taking the average over $i=1,\dots,N$ we get the thesis.
\qed\end{proof}
From now on, we denote with $D_0$ the following quantity:
\be \label{D0}
D_0 = \sum_{s=1}^d\E\lrq{\diam(\ybf^s_0)^2}.
\ee
Lemma \ref{lemmaAk} and Lemma \ref{lemmaBk} below provide an upper bound to the expectation of $A_k$ and $B_k$. The proof of Lemma \ref{lemmaAk} is part of the proof of Theorem 3.2 in \cite{discreteCBO}. We include it in the appendix for completeness.
\begin{lemma}\label{lemmaAk}
For every $k\in\N$ let $A_k$ be the quantity defined in Lemma \ref{lemma14}. Assume that Assumptions A1 and A2 hold and $\ca<1$. Then 
\be\label{expAk_final}
    \E[A_k] \leq \Gamma_A e^{-\alpha f_*} e^{-2k(1-\ca)},
    \ee
    with $\Gamma_A =  M_H (1+(1-\gamma)^2+\xi^2)^{1/2}(\gamma^2+\xi^2)^{1/2} D_0$ where $D_0$ is given in \eqref{D0}.
\end{lemma}

\begin{lemma}\label{lemmaBk}
For every $k\in\N$ let $B_k$ be the quantity defined in Lemma \ref{lemma14}. Assume that Assumptions A1 and  A5 hold and $\ca<1$. Then 

\be\label{expBk_final}
\E[B_k]\leq \Gamma_B e^{-\alpha f_*} e^{-k(1-\ca)} 
\ee
with
$\Gamma_B=
M_g(\gamma\alpha M_v + \xi)  D_0^{1/2}$,
$M_g$ given in Corollary \ref{boundedfgf} and $D_0$  in \eqref{D0}.  

\end{lemma}
\begin{proof}
    For every $s=1,\dots,d$ we have
$$
\sum_{i=1}^Ne^{-\alpha f(\xbf^i_k)}(x^i_{k,s} -x^\alpha_{k,s}) = \sum_{i=1}^Ne^{-\alpha f(\xbf^i_k)}\lr{x^i_{k,s} - \frac{\sum_{j=1}^N{e^{-\alpha f(\xbf^j_k)}x^j_{k,s}}}{\sum_{j=1}^Ne^{-\alpha f(\xbf^j_k)}}}=0.
$$
Using this equality and \eqref{CBOhat_compatto} we have 
\bes
\begin{aligned}
B_k& = \frac{1}{N}\sum_{i=1}^Ne^{-\alpha f(\xbf^i_{k})}\nabla f(\hat\xbf^\alpha_k)^\top(\Delta \xbf)^i_k 
= -\frac{1}{N}\nabla f(\hat\xbf^\alpha_k)^\top\sum_{i=1}^Ne^{-\alpha f(\xbf^i_{k})}(\gamma I_d+\diageta_{k}^i)(x^i_{k}-\hat x^\alpha_{k})  \\   
& = -\frac{\nabla f(\hat\xbf^\alpha_k)^\top}{N}\lr{\gamma\sum_{i=1}^Ne^{-\alpha f(\xbf^i_{k})}\lr{(x^i_{k}-x^\alpha_{k})+(x^\alpha_{k}-\hat x^\alpha_{k})} 
+\sum_{i=1}^Ne^{-\alpha f(\xbf^i_{k})}\diageta_{k}^i(x^i_{k}-\hat x^\alpha_{k})} \\
& =- \frac{\nabla f(\hat\xbf^\alpha_k)^\top}{N}\lr{\gamma \sum_{i=1}^Ne^{-\alpha f(\xbf^i_{k})} (x^\alpha_{k}-\hat x^\alpha_{k}) 
+\sum_{i=1}^Ne^{-\alpha f(\xbf^i_{k})}\sum_{s=1}^d\eta_{k,s}^i(x^i_{k,s}-\hat x^\alpha_{k,s})\ebf_s} \\
&\leq e^{-\alpha f_*}\|\nabla f(\hat\xbf^\alpha_k)\|\lr{\gamma\|x^\alpha_{k} - \hat x^\alpha_{k}\|
 + \frac{1}{N}\sum_{i=1}^N\lr{\sum_{s=1}^d(\eta^i_{k,s})^2(x^i_{k,s} - \hat x^\alpha_{k,s})^2}^{1/2}}.
\end{aligned}
\ees
Taking the expectation in the inequality above, by Holder's inequality we have
$$
\begin{aligned}
\E[B_k]&\leq \gamma e^{-\alpha f_*}\E[\|\nabla f(\hat\xbf^\alpha_k)\|^2]^{1/2}\E\lrq{\|x^\alpha_{k} - \hat x^\alpha_{k}\|^2}^{1/2} \\
&\quad + \frac{e^{-\alpha f_*}}{N}\sum_{i=1}^N\E[\|\nabla f(\hat\xbf^\alpha_k)\|^2]^{1/2}\E\lrq{\sum_{s=1}^d(\eta^i_{k,s})^2(x^i_{k,s} - \hat x^\alpha_{k,s})^2}^{1/2}\\
&\leq  e^{-\alpha f_*}M_g\lr {\gamma \E\lrq{\|x^\alpha_{k} - \hat x^\alpha_{k}\|^2}^{1/2} 
+ \frac{\xi}{N}\sum_{i=1}^N\E\lrq{\|x^i_{k} - \hat x^\alpha_{k}\|^2}^{1/2}},\\
\end{aligned}
$$
where in the last inequality we used Corollary \ref{boundedfgf} $iii)$.

By Lemma \ref{lemmadiffxalpha} and \ref{exp_max_err}, inequalities  \eqref{diam_rec} and \eqref{expmax}   we have for  $s=1,\dots,d$:
\be\label{derivazione_va}
\begin{aligned}
&\E\lrq{(x^\alpha_{k,s}-\hat x^\alpha_{k,s})^2}    \leq \alpha^2\E\lrq{\|{\cal E}_k\|_{\infty}^2\diam(\ybf^s_0)^2\prod_{l=0}^{k-1}(1-\gamma+4\max_{i=1:N}|\eta_{l,s}^i|)^2}   \\  &\quad \quad \quad= \alpha^2\E\lrq{\E\lrq{\|{\cal E}_k\|_{\infty}^2\diam(\ybf^s_0)^2\prod_{l=0}^{k-1}(1-\gamma+4\max_{i=1:N}|\eta_{l,s}^i|)^2\big |{\sigma_{k}} }}   \\  
& \quad \quad \quad\leq \alpha^2{M_v^2}\E\lrq{\diam(\ybf^s_0)^2\prod_{l=0}^k(1-\gamma+4\max_{i=1:N}|\eta_{l,s}^i|)^2}\\  & \quad \quad \quad\leq \alpha^2{M_v^2} \ca^{2k} \E\lrq{\diam(\ybf^s_0)^2} \leq \alpha^2{M_v^2} e^{-2k(1-\ca)}\E\lrq{\diam(\ybf^s_0)^2},
\end{aligned}
\ee
where in the last inequality we used that $e^t-1\geq t$ for every $t\in\R$.
Therefore 
\begin{equation}\label{expaux}
    \E\lrq{\|x^\alpha_{k} - \hat x^\alpha_{k}\|^2}^{1/2}  
\leq \alpha {M_v} e^{-k(1-\ca)} D_0^{1/2}
\end{equation}
 and part $iv)$ of Lemma \ref{lemma_diam_CBO} yields
the thesis.
\qed\end{proof}

We are now in the position to prove the main result of this subsection that provides the sought upper bound on the optimality gap. 
The result is an extension to the noisy case to Theorem 3.2 in \cite{discreteCBO}. The main difference is in  the second term of the right hand side of \eqref{ass_initialdistr}, which is necessary to take into account the noise in the objective function evaluation.
We will make use of the following result given in \cite{discreteCBO}.
\begin{lemma}\label{lemmaf}\cite[Proposition 3.1]{discreteCBO} 
    Let  Assumptions A1, A3 and A4 hold, and let $\xbf_0$ denote a random variable with law $\rho_0$ given in Assumption 4. Then, for $\alpha$ sufficiently large,
    $$-\frac{1}{\alpha}\log\lr{\E\lrq{e^{-\alpha f(\xbf_0)}}} = f_{*}+\frac{d}{2}\frac{\log(\alpha)}{\alpha} +\mathcal{O}\lr{\frac{1}{\alpha}}.$$
\end{lemma}

\begin{theorem}\label{thmessinf}
    Let us assume that $A1-A5$ hold,and let $\xbf_0$ denote a random variable with law $\rho_0$ given in Assumption A4. 
    Moreover, let us assume that the parameters $\gamma, \xi$ of the methods yield $\ca<1$, and that the initial distribution is such that for some $\varepsilon>0$
    \be\label{ass_initialdistr}
    \begin{aligned}
 (1-\varepsilon)\E\lrq{e^{-\alpha f(\xbf_0)}}&\geq \frac{{\alpha e^{-\alpha f_*}}\Gamma_A}{1-e^{-2(1-\ca)}} 
    +\frac{{\alpha e^{-\alpha f_*}}\Gamma_B }{1-e^{-(1-\ca)}},
    \end{aligned}
    \ee
where $\Gamma_A$ and $\Gamma_B$ are given in Lemma \ref{lemmaAk} and Lemma \ref{lemmaBk}, respectively.
    Then, given the point $\xbf_{\infty}$ defined in Theorem \ref{convergencethm}, it holds
    \begin{equation}\label{bound_essinf}\left|\essinf f(\xbf_\infty)-f_*\right|<\frac{d}{2}\frac{\log(\alpha)}{\alpha}+ G(\alpha) \end{equation}
    for some function $G(\alpha)=
\mathcal{O}\lr{\frac{1}{\alpha}}$.
\end{theorem}
\begin{proof}
 Taking the expectation in \eqref{meandiffexp}, using Lemma \ref{lemmaAk} and \ref{lemmaBk} we have
   
       $$
       \begin{aligned}
    \frac{1}{N}\sum_{i=1}^N\E\lrq{e^{-\alpha f(\xbf^i_{k+1})}}
    &\geq \frac{1}{N}\sum_{i=1}^N \E\lrq{e^{-\alpha f(\xbf^i_{k})}} -\alpha \E[A_k]-\alpha \E[B_k]  \\  
    &\geq \frac{1}{N}  \sum_{i=1}^N \E\lrq{e^{-\alpha f(\xbf^i_{k})}}     
    -\alpha e^{-\alpha f_*}\left(\Gamma_A\ \zeta_k^2 
   + \Gamma_B \zeta_k\right)
        \end{aligned}
    $$
    where $\zeta_k=e^{-k(1-\ca)}$.
Recursively applying this inequality we get
       \be\label{avgexp}
       \begin{aligned}
    \frac{1}{N}\sum_{i=1}^N\E\lrq{e^{-\alpha f(\xbf^i_{k+1})}}
    &\geq \frac{1}{N}\sum_{i=1}^N \E\lrq{e^{-\alpha f(\xbf^i_{0})}}   
    -\alpha e^{-\alpha f_*}\left(\Gamma_A \sum_{l=0}^k \zeta_l^2+
  \Gamma_B\sum_{l=0}^k \zeta_l\right).
        \end{aligned}
    \ee
    Since all the particles $\xbf^i_0$ are random variables with law $\rho_0$, which is also the law of the random variable $\xbf_0$, we have 
    \begin{equation*}
        \frac{1}{N}\sum_{i=1}^N \E\lrq{e^{-\alpha f(\xbf^i_{0})}} = \frac{1}{N}\sum_{i=1}^N \E\lrq{e^{-\alpha f(\xbf_{0})}} =\E\lrq{e^{-\alpha f(\xbf_{0})}}. 
    \end{equation*}
    Since $\ca<1$ by hypothesis, we have
$$    \sum_{l=0}^\infty \zeta_l = \frac{1}{1-e^{-(1-\ca)}}=\frac{1}{1-\zeta_1},\quad\quad
    \sum_{l=0}^\infty \zeta_l^2 = \frac{1}{1-e^{-2(1-\ca)}}=\frac{1}{1-\zeta_1^2},
  $$
therefore, taking the limit as $k\rightarrow+\infty$ in \eqref{avgexp}
     \be
       \begin{aligned}
    \E\lrq{e^{-\alpha f(\xbf_\infty)}} & = \frac{1}{N}\sum_{i=1}^N\E\lrq{e^{-\alpha f(\xbf_\infty)}} =\lim_{k\rightarrow+\infty} \frac{1}{N}\sum_{i=1}^N\E\lrq{e^{-\alpha f(\xbf_{k+1}^i)}} \\
    & \geq \E\lrq{e^{-\alpha f(\xbf_{0})}}  
    -\alpha e^{-\alpha f_*}\left(\frac{\Gamma_A} {1-\zeta_1^2}+\frac{\Gamma_B }{1-\zeta_1}\right).
        \end{aligned}
    \ee
Applying  \eqref{ass_initialdistr} we conclude that 
    \be
    \E\lrq{e^{-\alpha f(\xbf_\infty)}}
    \geq \varepsilon\E\lrq{e^{-\alpha f(\xbf_{0})}}
    \ee
    By definition of essential limit we then have 
    $$e^{-\alpha\essinf f(\xbf_\infty)} = \E\lrq{e^{-\alpha\essinf f(\xbf_\infty)}}\geq\E\lrq{e^{-\alpha f(\xbf_\infty)}}\geq \varepsilon \E\lrq{e^{-\alpha f(\xbf_{0})}}. $$
    Taking the logarithm on both sides of the previous inequality and applying Lemma \ref{lemmaf} we have 
    \be
       \begin{aligned}
       \essinf f(\xbf_\infty)&\leq -\frac{\log(\varepsilon)}{\alpha}-\frac{1}{\alpha}\log\lr{\E\lrq{e^{-\alpha f(\xbf_{0})}}}  \\  
       & <f_*+\frac{d\log(\alpha)}{2\alpha}+G(\alpha) 
        \end{aligned}
        \ee
     for some function $G(\alpha)=
\mathcal{O}\lr{\frac{1}{\alpha}}$,  which gives the thesis.
\qed\end{proof}

\begin{remark}
Let us consider condition \eqref{ass_initialdistr} in the theorem above. The left hand side is larger whenever $f(\xbf_0)$ is closer to $f_*$, while the right hand side is smaller if the initial points are closer to their average and if the error ${\cal E}_k$ in the function evaluation is smaller. As $\alpha$ increases, the condition becomes more restrictive, but the thesis ensures that $\essinf f(\xbf_\infty)$ is closer to $f_*$. 
\end{remark}
\begin{remark}
    Let us denote with $\widehat G(\alpha)$ the quantity $\frac{d}{2}\frac{\log(\alpha)}{\alpha}+ G(\alpha)$.
    We want to prove that Theorem \ref{thmessinf} implies 
    $$\mathbb{P}\left\{f(x_{\infty})-f_*<\widehat G(\alpha)\right\}>0.$$
    First of all, let us notice that inequality \eqref{bound_essinf} is equivalent to 
    $$\sup \big\{a\in\R\ |\ \essinf(f(x_{\infty})-f_*)= a\big\}< \widehat G(\alpha)$$
    which in turn is equivalent to
    \begin{equation}\label{auxineq}\sup \big\{a\in\R\ |\ \mathbb{P}\{(f(x_{\infty})-f_*)\geq a\}=1\big\}< \widehat G(\alpha).\end{equation}
    Assume by contradiction that $\mathbb{P}\left\{f(x_{\infty})-f_*\geq\widehat G(\alpha)\right\}=1$. This implies
    $$\widehat G(\alpha)\in\big\{a\in\R\ |\ \mathbb{P}\{(f(x_{\infty})-f_*)\geq a\}=1\big\}$$
    and therefore
    $$\sup \big\{a\in\R\ |\ \mathbb{P}\{(f(x_{\infty})-f_*)\geq a\}=1\big\}\geq \widehat G(\alpha),$$
    which contradicts \ref{auxineq}.\\
    Let us assume the algorithm is run $R$ times, with the same parameters and the same initial distribution of the particles, and for $r=1,\dots,R$ let us denote with $S_r$ the random variable
    $$S_r = \begin{cases}1 & \text{if } f(x_{\infty})-f_*<\widehat G(\alpha) \\
    0 & \text{otherwise}\end{cases},$$
    and with $p$ the probability
    $p = \mathbb{P}\left\{S_r=1\right\},$ which we proved is strictly positive. We then have 
    $$\mathbb{P}\left\{\exists\ r\in\{1,\dots, R\}\ | \ S_r=1\right\} = 1 - \mathbb{P}\left\{ S_r=0\ \forall k=1,\dots,R\right\} = 1-(1-p)^R,$$
    which tends to 1 as $R$ tends to infinity. Then, increasing the number of times
the algorithm is run, the probability of obtaining an accurate approximation of $f_*$ increases. However we notice that from Theorem 3.3 one cannot derive a lower bound for the probability $p$ of obtaining an accurate approximation of $f_*.$ Moreover, these inequalities concern the function values $f(x_{\infty})$ and $f_*$ and they do not imply in general that $x_{\infty}$ will be close to $x_*$.
\end{remark}

\subsection{Convergence in expectation}
In the following we show that, provided the objective function satisfies some additional regularity assumptions, one can prove that the mean squared distance of the particles from the solution $\xbf^*$, i.e.
\be \label{Vk}
V_{k} = \frac{1}{N}\sum_{i=1}^N\|\xbf^i_k-\xbf^*\|^2\ee
goes to zero in expectation, without any assumptions on the initial distribution of the particles. 
We first prove a crucial inequality linking the expected value of the mean squared distance at two consecutive iterations. 

\begin{lemma}\label{lemmaVk}
    For every $k\in\N$ we have
        \be
       \begin{aligned}
       \E[V_{k+1}]&\leq(1-2\gamma +\gamma^2+\xi^2)\E[V_{k}]+ (\gamma^2+\xi^2)\E[\|\hat\xbf^\alpha_k-\xbf^*\|^2]  \\  
       &+2(\gamma^2+\xi^2+(\gamma^2+\xi^2)^{1/2})\E[V_{k}]^{1/2}\E[\|\hat\xbf^\alpha_k-\xbf^*\|^2]^{1/2}.
       \end{aligned}
    \ee
\end{lemma}
\begin{proof}
    By definition of $V_{k+1}$ and \eqref{CBOhat_compatto} we have
     \be\label{Vk+1}
       \begin{aligned}
       V_{k+1}& = \frac{1}{N}\sum_{i=1}^N\|\xbf^i_{k+1}-\xbf^*\|^2=  
       \frac{1}{N}\sum_{i=1}^N\left\|\xbf^i_k+(\Delta x)_k^i-\xbf^*\right\|^2   \\  
       &= \frac{1}{N}\sum_{i=1}^N\|\xbf^i_k-\xbf^*\|^2+\frac{1}{N}\sum_{i=1}^N\|(\Delta x)_k^i\|^2 +\frac{2}{N}\sum_{i=1}^N(\xbf^i_k-\xbf^*)^\top (\Delta x)_k^i\\
       &= V_k+Z_{k,1}+Z_{k,2},
       \end{aligned}
    \ee
with
    \be
       \begin{aligned}
       &Z_{k,1} = \frac{1}{N}\sum_{i=1}^N\|(\Delta x)_k^i\|^2= \frac{1}{N}\sum_{i=1}^N\left\|(\gamma I_d+\diageta_{k}^i)(x^i_{k}-\hat x^\alpha_{k})\right\|^2,\\
       &Z_{k,2} =  \frac{2}{N}\sum_{i=1}^N(\xbf^i_k-\xbf^*)^\top (\Delta x)_k^i\\
       &\quad \quad =-\frac{2}{N}\sum_{i=1}^N(\xbf^i_k-\xbf^*)^\top(\gamma I_d+\diageta_{k}^i)(x^i_{k}-\hat x^\alpha_{k}).
       \end{aligned}
    \ee
    By summing and subtracting  $x^*$ and using Cauchy-Schwartz inequality we get
    \be\label{Ak2}
       \begin{aligned}
       Z_{k,1}
        &\leq \frac{1}{N}\sum_{i=1}^N\left\|(\gamma I_d+\diageta_{k}^i)(x^i_{k}-\xbf^*)\right\|^2  \\  
       &+\frac{2}{N}\sum_{i=1}^N\left\|(\gamma I_d+\diageta_{k}^i)(x^i_{k}-\xbf^*)\right\|\left\|(\gamma I_d+\diageta_{k}^i)(\hat x^\alpha_{k}-\xbf^*)\right\|  \\  
       &+\frac{1}{N}\sum_{i=1}^N\left\|(\gamma I_d+\diageta_{k}^i)(\hat x^\alpha_{k}-\xbf^*)\right\|^2.
       \end{aligned}
    \ee
    By Holder's inequality we have
    \bes
       \begin{aligned}
    &\E [
    \left\|(\gamma I_d+\diageta_{k}^i)(x^i_{k}-\xbf^*)\right\|\left\|(\gamma I_d+\diageta_{k}^i)(\hat x^\alpha_{k,}-\xbf^*)\right\|] 
     \\  
       & \le \lr{\sum_{s=1}^d\E[(\gamma+\eta_{k,s}^i)^2]\E[(x^i_{k,s}-\xbf^*_s)^2]}^{1/2}\lr{\sum_{s=1}^d\E[(\gamma+\eta_{k,s}^i)^2]\E[(\hat x^\alpha_{k,s}-\xbf^*_s)^2]}^{1/2}  \\  
       &=(\gamma^2+\xi^2)\lr{\sum_{s=1}^d\E[(x^i_{k,s}-\xbf^*_s)^2]}^{1/2}\lr{\sum_{s=1}^d\E[(\hat x^\alpha_{k,s}-\xbf^*_s)^2]}^{1/2}  \\  
       &= (\gamma^2+\xi^2)\E[\|x^i_{k}-\xbf^*\|^2]^{1/2}\E[\|\hat x^\alpha_{k}-\xbf^*\|^2]^{1/2}.
       \end{aligned}
    \ees
    In addition we note that, given $\{a_i\}_{i=1}^N$ with $a_i\ge 0$ for any $i$, it holds 
        \be \label{suma}
        \frac{1}{N}\sum_{i=1}^N a_i\le \frac{1}{\sqrt{N}}\lr{\sum_{i=1}^N a_i^2}^{1/2},
        \ee
        that yields
        \be \label{suma1}
        \frac{1}{N}\sum_{i=1}^N\E[\|x^i_{k}- x^*\|^2]^{1/2}\le  \frac{1}{\sqrt{N}}\left(\sum_{i=1}^N\E[\|x^i_{k}- x^*\|^2]\right)^{1/2}=\E[V_k]^{1/2}.
        \ee
       Taking the expected value on both sides of \eqref{Ak2} and using \eqref{suma1} we then have
    \be\label{expAk2}
       \begin{aligned} 
       \E[Z_{k,1}]&  \le  \frac{1}{N}\sum_{i=1}^N\sum_{s=1}^d\E[(\gamma+\eta_{k,s}^i)^2]\E[(x^i_{k,s}-\xbf^*_s)^2]  \\  
       &\quad +\frac{2}{N}(\gamma^2+\xi^2)\sum_{i=1}^N\E[\|x^i_{k}-\xbf^*\|^2]^{1/2}\E[\|\hat x^\alpha_{k}-\xbf^*\|^2]^{1/2}\\
       &\quad +\frac{1}{N}\sum_{i=1}^N\sum_{s=1}^d\E[(\gamma+\eta_{k,s}^i)^2]\E[(\hat x^\alpha_{k,s}-\xbf^*_s)^2]   \\  
       &\le (\gamma^2+\xi^2)(\E[V_k]+2\E[V_k]^{1/2}\E[\|\hat x^\alpha_{k}-\xbf^*\|^2]^{1/2}\\
       &\quad +\E[\|\hat x^\alpha_{k}-\xbf^*\|^2]).\\
       \end{aligned}
         \ee
    We now consider the term $Z_{k,2}$ and again   sum and subtract  $x^*$: 
    \be
       \begin{aligned}
       Z_{k,2} 
     &\leq -\frac{2}{N}\sum_{i=1}^N\sum_{s=1}^d(\gamma+\eta_{k,s}^i)(x^i_{k,s}- x^*_{s})^2  \\  
       &\quad +\frac{2}{N}\sum_{i=1}^N\|\xbf^i_k-\xbf^*\|\left\|(\gamma I_d+\diageta_{k}^i)(\hat x^\alpha_{k}-x^*)\right\|.\\
       \end{aligned}
    \ee
    Using Holder's inequality and proceeding as in \eqref{expAk2} we get
    \be \label{expBk2}
       \begin{aligned}
       \E[Z_{k,2}] \le & -2\gamma\frac{1}{N}\sum_{i=1}^N\E[\|x^i_{k}- x^*\|^2]  \\ &+2(\gamma^2+\xi^2)^{1/2}\E[\|\hat\xbf^\alpha_k-\xbf^*\|^2]^{1/2}\frac{1}{N}\sum_{i=1}^N\E[\|x^i_{k}- x^*\|^2]^{1/2} \\  
       \le &-2\gamma\E[V_k]+2(\gamma^2+\xi^2)^{1/2}\E[\|\hat\xbf^\alpha_k-\xbf^*\|^2]^{1/2}\E[V_k]^{1/2}.
       \end{aligned}
       \ee
    From \eqref{Vk+1}, \eqref{expAk2} and \eqref{expBk2} we get the thesis.
\qed\end{proof}
Given a scalar $r>0$ we now introduce the ball of radius $r$ centered in $\xbf^*$: 
\be\label{Br}
\B_r:=\{\xbf\in\R^n\enspace|\enspace \|\xbf-\xbf^*\|\leq r \}
\ee
and the following set of indices
\be \label{Ir}
{\cal I}_{k, r}=\{i \in \{1\ldots,N\}\, |\, \xbf_k^i \in \B_{ r}\}.
\ee
Inspired by the analysis in the infinite dimensional setting carried out in \cite{Fornasier1},
we now make an additional assumption on the objective function $f.$

\noindent\textbf{Assumption A6.} There exist constants $f_{\infty}, R_0, \beta, \nu>0$ such that
\begin{equation*}\begin{aligned}
    &(f(\xbf)-f_*)^{\nu}\geq \beta\|\xbf-\xbf^*\|\enspace \text{for every}\ \xbf\in \B_{R_0}\\
    &f(\xbf)\geq f_*+f_\infty\enspace \text{for every}\ \xbf\in \B_{R_0}^\complement.
\end{aligned}
\end{equation*}

\begin{remark}
    Assumption A6 poses condition on how well isolated the global minimizer is, in addition to the uniqueness assumption A3. The first inequality implies that in a ball around $x_*$, the function cannot grow too slowly. Namely, the difference $f(x)-f_*$ has to grow at least as $(\beta\|x-x_*\|)^{1/\nu}$. The second inequality ensures that, outside of a neighborhood of $x_*$, the function cannot take values that are too close to the minimum $f_*.$
\end{remark}

The following Lemma, known in the literature as \emph{quantitative Laplace principle}, provides
an upper bound on $\|\xbf_k^\alpha-\xbf^*\|$, where $\xbf_k^\alpha$ has been defined in \eqref{xkalpha}. 
Consider $r\in(0,R_0]$ and $q>0$ such that $q+f_r\leq f_{\infty}$, with 
\be \label{fr}
f_r = \max_{x \in \B_r} f(x).
\ee
We distinguish the case 
$|{\cal I}_{k,r}|>0$ and $|{\cal I}_{k,r}|=0$. We stress that in the first case, at iteration $k$  there exists at least one particle that belongs to $B_r$. 
Part (a) of the Lemma is completely analogous to Proposition 1 in \cite{Fornasier1}, while part (b) is obtained from (a) with minor changes. For completeness, a proof is provided in the appendix.
\begin{lemma}[Quantitative Laplace Principle]\label{diffxalphasol}
Let Assumptions $A1-A6$ hold and $x_k^\alpha$ given in   \eqref{xkalpha}. Consider $r\in(0,R_0]$ and $q>0$ such that $q+f_r\leq f_{\infty}$, with 
$f_r$ given in \eqref{fr} and 
 $\B_r $ in \eqref{Br}.
For every $k\in\N$ we have
\begin{enumerate}[a)]
    \item if  $|{\cal I}_{k,r}|>0$ then
\begin{equation*}
    \|\xbf^\alpha_k-\xbf^*\|\leq \frac{1}{\beta}(q+f_r)^\nu+\frac{e^{-\alpha q}}{ |{\cal I}_{k,r}|}\sum_{i=1}^N\|\xbf^i_k-\xbf^*\|,
\end{equation*}
    \item if  $|{\cal I}_{k,r}|=0$ it still holds that
    $$\|\xbf^\alpha_k-\xbf^*\|\leq \frac{1}{\beta}(q+f_r)^\nu+\frac{e^{-\alpha q}}{e^{-\alpha\max_{i=1:N} f(\xbf^i_k)}}\frac{1}{N}\sum_{i=1}^N\|\xbf^i_k-\xbf^*\|.$$
\end{enumerate}

\end{lemma}
In the following theorem, given $\varepsilon>0$  and assuming that at a specific iteration  $\bar k$  the condition  $|{\cal I}_{\bar k,r}|>0 $ holds,
we provide an upper bound to the number of iterations needed to reach $\E[V_{k}]\leq \varepsilon$,
where $V_k$ is the mean squared distance from $x^*$ given in \eqref{Vk}.
Notice that $\ca<1$ and $\gamma \in (0,1]$ imply $1-(1-\gamma)^2-\xi^2>0$.

\begin{theorem}\label{thmVkdet}
Le us assume that the assumptions in Lemma \ref{diffxalphasol} hold and that  $\ca<1$. For any $\varepsilon>0,\ \tau\in(0,1)$ let us define
 \be\label{constants}
       \begin{aligned}
       &\mu = (1-(1-\tau)(1-(1-\gamma)^2-\xi^2))\in(0,1)\\
       &\hat\sigma  = \min\left\{\frac{\tau}{4}\frac{1-(1-\gamma)^2-\xi^2}{\gamma^2+\xi^2+(\gamma^2+\xi^2)^{1/2}}, \lr{\frac{\tau}{2}\frac{1-(1-\gamma)^2-\xi^2}{\gamma^2+\xi^2}}^{1/2}\right\}\\
       & C_k = \E[V_k]^{1/2}\hat\sigma,
       \end{aligned}
    \ee
where $V_k$ has been defined in \eqref{Vk}.
Let $\bar k$ such that
\be \label{bark}
\E[V_k]>\varepsilon,\quad\quad \E[\|\hat \xbf^\alpha_k-\xbf^*\|^2]\leq C_{k}^2\quad \forall k=0,\dots,\bar k -1
\ee
and 
\be\label{kstar}
k_* = \left\lceil \log_{1/\mu}\lr{\frac{\E[V_0]}{\varepsilon}}\right\rceil.
\ee
Moreover, assume that:
\be \label{rho_bark}
|{\cal I}_{\bar k,r}|>0.
\ee
Then,
for $\alpha$ sufficiently large and $M_v$ given in Lemma \ref{exp_max_err},
sufficiently small, it holds 
$$\E[V_{k_*}]\leq \varepsilon.$$ 
In addition, $\E[V_{k}]\leq\mu^k\E[V_0]$ for $k=0,\dots,k_*.$ 
\end{theorem}

\begin{proof}
By definition of $\bar k$ in \eqref{bark}, Lemma \ref{lemmaVk}, and the definition of $C_k$ and  $\hat\sigma$ in \eqref{constants} we have that for every $k=0,\dots,\bar k-1$
    \be\label{Vkmu}
       \begin{aligned}
       \E[V_{k+1}]\leq&(1-2\gamma +\gamma^2+\xi^2)\E[V_{k}]+ (\gamma^2+\xi^2)C_k^2  \\  
       &+2(\gamma^2+\xi^2+(\gamma^2+\xi^2)^{1/2})\E[V_{k}]^{1/2}C_k  \\  
       \leq&(1-2\gamma +\gamma^2+\xi^2)\E[V_{k}]+ \hat\sigma^2(\gamma^2+\xi^2)\E[V_{k}]  \\  
       &+2\hat\sigma(\gamma^2+\xi^2+(\gamma^2+\xi^2)^{1/2})\E[V_{k}]  \\  
       \leq&(1-(1-\tau)(1-(1-\gamma)^2-\xi^2))\E[V_k] \\
       =&\mu\E[V_k]\leq \mu^{k+1}\E[V_0].
       \end{aligned}
    \ee
If $k_*<\bar k$ then for every $k=0,\dots,k_*$ we have $\E[V_k]\leq \mu^{k+1}\E[V_0]$. In particular, by definition of $k_*$, $\E[V_{k_*}]\leq \varepsilon,$ and the thesis holds.  
Note that $\bar k$ is the first iteration where either $\E[V_{\bar k}]\leq\varepsilon$ or $\E[\|\hat\xbf^\alpha_{\bar k}-\xbf^*\|^2]>C_{\bar k}^2$.
If $k_*\geq \bar k$ and $\E[V_{\bar k}]\leq\varepsilon$, then there is nothing to prove. In case   $\E[V_{\bar k}]>\varepsilon$
we will show that $\E[\|\hat\xbf^\alpha_{\bar k}-\xbf^*\|^2]>C_{\bar k}^2$ cannot happen.
 Let us  define 
    \be \label{q_r}
       \begin{aligned}
       &q = \frac{1}{2}\min\left\{f_{\infty},\ \lr{\frac{\beta}{2}C_{\bar k}}^{1/\nu}\right\}\\
       &r = \max\{s\in[0,R_0]\enspace | \enspace f(\xbf)\leq q\ \forall \xbf\in B_s\}.
       \end{aligned}
    \ee

  {By definition we have that $r\leq R_0$, $q+f_r\leq f_\infty,$ and $q+f_r\leq 2q$.} Moreover, since $q>0$ and $f$ is continuous, $r$ is strictly positive.
    By Lemma \ref{diffxalphasol}, \eqref{suma}, \eqref{q_r} and \eqref{rho_bark} we have
    \be\label{diff1}
    \|\xbf^\alpha_{\bar k}-\xbf^*\|\leq \frac{1}{\beta}(q+f_r)^\nu+\frac{e^{-\alpha q}}{|{\cal I}_{\bar k,r}|}\sum_{i=1}^N\|\xbf^i_{\bar k}-\xbf^*\|  
    \leq \frac{1}{2}C_{\bar k}+Ne^{-\alpha q}V_{\bar k}^{1/2}.
    \ee
    Lemma \ref{lemmadiffxalpha} and \eqref{diam_rec} yield
  \be\label{diff2}
       \begin{aligned}
    \|\xbf^\alpha_{\bar k}-\hat\xbf^\alpha_{\bar k}\|&\leq \alpha \|{\cal E}_{\bar k}\|_{\infty}\lr{\sum_{s=1}^d\diam(\ybf^s_{\bar k})^2}^{1/2} \\
    &\leq\alpha\|{\cal E}_{\bar k}\|_{\infty} \lr{\sum_{s=1}^d\diam(\ybf^s_0)^2\prod_{l=0}^{{\bar k}-1}(1-\gamma+4\max_{i=1:N}|\eta_{l,s}^i|)^2}^{1/2}.
       \end{aligned}
    \ee
  
    Then,
    \be\label{eq:deriv}
       \begin{aligned}
\|\hat\xbf^\alpha_{\bar k}-\xbf^*\|^2 \le&
    \|\xbf^\alpha_{\bar k}-\xbf^*\|^2+\|\xbf^\alpha_{\bar k}-\hat\xbf^\alpha_{\bar k}\|^2 + 2\|\xbf^\alpha_{\bar k}-\xbf^*\|\|\xbf^\alpha_{\bar k}-\hat\xbf^\alpha_{\bar k}\|  \\  
   \leq& \frac{1}{4}C_{\bar k}^2+N^2e^{-2\alpha q}V_{\bar k}+ Ne^{-\alpha q}V_{\bar k}^{1/2}C_{\bar k} + \\  
    & \alpha^2 \|{\cal E}_{\bar k}\|_{\infty}^2 \sum_{s=1}^d\diam(\ybf^s_0)^2\prod_{l=0}^{{\bar k}-1}(1-\gamma+4\max_{i=1:N}|\eta_{l,s}^i|)^2+ \\  
    &C_{\bar k}\alpha\|{\cal E}_{\bar k}\|_{\infty} \lr{\sum_{s=1}^d\diam(\ybf^s_0)^2\prod_{l=0}^{{\bar k}-1}(1-\gamma+4\max_{i=1:N}|\eta_{l,s}^i|)^2}^{1/2} +  \\  
    &2\alpha Ne^{-\alpha q}V_{\bar k}^{1/2} \|{\cal E}_{\bar k}\|_{\infty} \lr{\sum_{s=1}^d\diam(\ybf^s_0)^2\prod_{l=0}^{{\bar k}-1}(1-\gamma+4\max_{i=1:N}|\eta_{l,s}^i|)^2}^{1/2}.
       \end{aligned}
    \ee
Taking the conditioned expected value, using Lemma \ref{exp_max_err}, \eqref{D0}, \eqref{prodexp}, the Holder's inequality, and proceeding as in \eqref{derivazione_va} 
    we obtain
 \begin{equation*}
           \begin{aligned}
               E[\|\hat\xbf^\alpha_{\bar k}-\xbf^*\|^2\enspace |\enspace |{\cal I}_{\bar k,r}|>0  ]&\leq 
    \frac{1}{4}C_{\bar k}^2+N^2e^{-2\alpha q}\E[V_{\bar k}]+ Ne^{-\alpha q}\E[V_{\bar k}]^{1/2}C_{\bar k}\\
    &\quad + \alpha^2{M_v^2} \ca^{2\bar k}D_0
     + \alpha {M_v}\ca^{\bar k}D_0^{1/2}
    \lr{C_{\bar k}+2Ne^{-\alpha q}\E\lrq{V_{\bar k}}^{1/2}}.\\
           \end{aligned}
       \end{equation*}
            Rearranging the terms 
    we get
     \be\label{expdist}
     \begin{aligned}
    &\E[\|\hat\xbf^\alpha_{\bar k}-\xbf^*\|^2 \enspace |\enspace |{\cal I}_{\bar k,r}|>0]  \leq 
    \frac{1}{4}C_{\bar k}^2+ C_{\bar k}
    \ca^{\bar k}D_0^{1/2}\alpha M_v + \ca^{2\bar k}D_0 \alpha^2M_v^2
      \\
      &\quad \quad \quad + e^{-\alpha q}N\E[V_{\bar k}]^{1/2}
       \lr{Ne^{-\alpha q}\E[V_{\bar k}]^{1/2} + C_{\bar k}+2 \ca^{\bar k}D_0^{1/2} \alpha M_v}.
       \end{aligned}
    \ee
Noticing that $C_{\bar k} = \hat\sigma\E[V_{\bar k}]^{1/2}\geq\hat\sigma\varepsilon^{1/2}$, it is easy to see that whenever
    \begin{equation}\label{Mvbound}
        M_v<\frac{\hat \sigma\varepsilon^{1/2} (\sqrt{3}-1)  }{2\alpha\ca^{\bar k}D_0^{1/2} }=:\frac{\widehat M}{\alpha}\end{equation}
        we have 
    \be \label{hatTK}
    C_{\bar k}
    \ca^{\bar k}D_0^{1/2} \alpha M_v +  \ca^{2\bar k}D_0 \alpha^2M_v^2\le \frac{1}{2}C_{\bar k}^2.
    \ee
    
Finally  we can find $\alpha_{min}>0$ so that  for any $\alpha\ge \alpha_{min}$
  \be
            e^{-\alpha q}N\E[V_{\bar k}]^{1/2}\lr{Ne^{-\alpha q}\E[V_{\bar k}]^{1/2} + C_{\bar k}+2N^{1/2}\ca^{\bar k}D_0^{1/2} \alpha\widehat M} \le 
           \frac{1}{4}C_{\bar k}^2.
          \ee
  Then, letting 
$\alpha=\hat \alpha>\alpha_{min}$
and $M_v\le \widehat M/\hat \alpha$, 
   it holds
    $$\E[\|\hat\xbf^\alpha_{\bar k}-\xbf^*\|^2]\leq C_{\bar k}^2,$$
    which concludes the proof.
\qed\end{proof}

 The same result as in Theorem \ref{thmVkdet} can be proved under weaker assumptions, i.e. assuming that condition \eqref{rho_bark} holds only with some fixed probability. The proof is included in the Appendix for completeness.
\begin{theorem}\label{thmVkprob}
Let us assume that $\ca<1$, the assumptions in Lemma \ref{diffxalphasol} hold and 
    \begin{equation}\label{rho_bark_prob}
        \mathbb P\left\{|{\cal I}_{\bar k,r}|>0\right\}\ge 1-\delta \enspace \text{with}\ \delta \leq\hat\delta e^{-4\alpha M_f} \end{equation}
    where $\bar k$ is given in \eqref{bark}, $M_f$ is the upper bound from Corollary \ref{boundedfgf}, and $\hat\delta\in (0,e^{4\alpha M_f})$. 
Then, for $\alpha$ sufficiently large and $M_v$ sufficiently small, it holds  $\E[V_{k_*}]\leq \varepsilon$, with $k_*$ given in \eqref{kstar}. Moreover, $\E[V_{k}]\leq\mu^k\E[V_0]$ for $k=0,\dots,k_*,$ with $\mu$ given in \eqref{constants}.
\end{theorem}

\begin{remark} Let us consider the condition \eqref{Mvbound} on $M_v$ from the previous theorems. 
In the general case, given the form of $M_V$ given in  Lemma \ref{exp_max_err}, we have that 
in order to satisfy 
\eqref{Mvbound}  it must hold
$$(t_0+t_1M_f^2)^{1/2}\leq\frac{\sigma\varepsilon^{1/2} (\sqrt{3}-1)  }{2\alpha\ca^{\bar k}D_0^{1/2}N^{1/2} }\sim\frac{1}{N^{1/2}}.$$
In case of additive Gaussian noise this condition becomes significantly less restrictive:
$$(t_0+t_1M_f^2)^{1/2}\leq\frac{\sigma\varepsilon^{1/2} (\sqrt{3}-1)  }{4\alpha\ca^{\bar k}D_0^{1/2}\log(\sqrt{2}N)^{1/2} }\sim\frac{1}{\log(\sqrt{2}N)^{1/2}}.$$
\end{remark}

    \begin{remark}
        From the theoretical analysis that we carried out, we can notice a trade-off between the exploration power of the system of particles and the convergence to a consensus point. Intuitively, if the number of particles $N$ is larger, they will visit more points in the search space, while it will be harder for the system to achieve consensus. This can be seen more rigorously in Theorems \ref{thm_cons1}, \ref{thmVkdet} and \ref{thmVkprob} as condition \ref{rho_bark} is more likely to hold if the number of particles $N$ is large while to ensure consensus one needs $\ca = 1-\gamma+8\xi(\log(\sqrt{2}N))^{1/2}<1$, which is more restrictive as $N$ increases. {
        In the numerical result section we will show that the condition $\ca<1$ appears in practice to be too restrictive.\\
        Assumption \eqref{rho_bark}, or its probabilistic relaxation \eqref{rho_bark_prob}, are necessary for our analysis, however in \cite{Fornasier1} it has been shown that for the continuous algorithm the analogous of $|{\cal I}_{k,r}|$, namely $\rho_t(\mathcal{B}_r)$, is guaranteed to be strictly positive for $t\in[0,T]$ for any $T$. }
    \end{remark}

\section{Numerical Results}\label{sec:numres}
In this section we investigate on the sensitivity of our approach to noise in the objective function evaluation. We first consider the minimization of the well known Rastrigin function with Gaussian noise on the objective function and then we analyze the behaviour of CBO applied to finite sum minimization problems employing  subsampling estimators of the objective function.

\subsection{Rastrigin Function}
We consider the global optimization problem $\min_{\R^d} f(x)$ with $f(x)=R(x)$, where $R$ is the Rastrigin function, defined as
\begin{equation}\label{Rastrigin}
    R(\xbf) = \sum_{i=1}^d\lr{x_i^2-10\cos(2\pi x_i)}+10^d,
\end{equation}
which is known to have global minimum at $\xbf^* = 0$ and several local minima, and focus on the computation of its global minimizer. Given a point $\xbf$ we define the oracle $\hat f$ by adding additive Gaussian noise of the form \eqref{Gaussian_noise},
where $\omega(\xbf)\in\R^2$ is a random variable with Gaussian components $\omega_j(\xbf)\sim\mathcal{N}(0,\sigma_j^2)$ for $j=0,1$.
When $\sigma_0=0$ we consider relative  random Gaussian noise, while the choice $\sigma_1=0$ produces absolute random Gaussian noise. We refer to the case where  both $\sigma_0$ and $\sigma_1$ are greater than zero as mixed Gaussian noise.

To study how the behaviour of the method changes with respect to the accuracy of the oracle and the parameter $\alpha$, we apply the CBO method \eqref{CBOhat} to the Rastrigin function with different values for the noises' standard deviations $\sigma_0,\ \sigma_1$ and different values of $\alpha\in\{10^4,10^3,10^2,50,10,5,1,0.1\}$. 
For every combination of the noises and $\alpha$ we run the method with parameters $\gamma = 0.1$, $\xi = 0.0056$, for 1000 iterations and 1000 runs. {Notice that this choice of the parameters satisfies the condition $\ca<1$ required for the convergence analysis presented in the previous sections.}
As a measure of accuracy, we consider the distance of the consensus point at the final iteration from the solution $\xbf^*=0.$

In Table \ref{table:Rastrigin1_new} we report our findings  for dimension $d=1$ and $N=100.$ 
For each couple $(\sigma_0, \sigma_1)$, we take the value of $\alpha$ that yields the best results.
We report the value of $\sigma_0$, $\sigma_1$, $\alpha$, the average, minimal and maximal distance from the solution (denoted with \emph{mean err.}, \emph{min err.} and \emph{max err.}, respectively), and the $50, 75$ and $90$-th percentile, i.e. the distance from the solution reached by at least $50\%$, $75\%$ and $90\%$ of the runs, respectively.

We can observe that the absolute Gaussian noise affects the quality of the solution as expected, but we manage to obtain an average error of the order of $10^{-3}$ (against $10^{-5}$ obtained in the noise-free case), even when an high noise level is present. We also underline that 
the method is quite insensitive to this type of noise as increasing $\sigma_0$ we observe only a slight decrease of the overall obtained accuracy.
Regarding the relative Gaussian noise, we note that the noise does not affect the accuracy of the computed approximation for values of the noise level smaller or equal than $0.25$,  while the method is not able to handle higher values of the noise. 
In case of the mixed noise we obviously observe the effect of both the absolute and relative noise. 
Overall the observed behaviour is consistent with our theoretical analysis.
Indeed, under Assumption A5, Theorems \ref{thmVkdet} and \ref{thmVkprob} shows that global convergence is obtained provided that the the upper bound   $t_0+t_1 M_f^2$ on the noise level is sufficiently small. In case of relative error, the bound on  $\sigma_1$ is more severe than the bound on $\sigma_0$ at particles $x$ at which $f$ takes large value. \\

\begin{table}
\begin{tabular}{ c |c| c| c| c| c| c| c| c }
$\sigma_0$ & $\sigma_1$ & $\alpha$ & mean err. & min err. & max err. &  50 perc. & 75 perc & 90 perc\\
\hline
\multicolumn{9}{c}{without noise  }  \\
\hline
0 &0   & $10^4$& 4.97e-5& 1.49e-8& 1.91e-4& 3.82e-5& 7.40e-5& 1.00e-4\\
\hline
\multicolumn{9}{c}{absolute noise  }  \\
\hline
0.1 &0 & 10& 1.62e-3& 7.99e-7& 6.98e-3& 1.32e-3& 2.41e-3& 3.41e-3\\
0.25 &0& 10& 2.21e-3& 1.42e-6& 9.15e-3& 1.83e-3& 3.09e-3& 4.67e-3\\
0.5 &0 & 5& 3.49e-3& 1.52e-5& 1.66e-2& 2.87e-3& 5.08e-3& 7.21e-3\\
1.0 &0 & 5& 7.25e-3& 1.06e-5& 3.48e-2& 6.00e-3& 1.04e-2& 1.47e-2\\
\hline
\multicolumn{9}{c}{relative noise  }  \\
\hline
0 &0.1 & $10^4$& 5.01e-5& 2.90e-8& 1.94e-4& 3.91e-5& 7.40e-5& 1.02e-4\\
0 &0.25& $10^4$& 5.45e-5& 7.80e-9& 2.55e-4& 4.81e-5& 7.69e-5& 1.15e-4\\
0 &0.5 & 0.1   & 1.14e-1& 3.29e-5& 7.50e-1& 7.81e-2& 1.49e-1& 2.53e-1\\
\hline
\multicolumn{9}{c}{mixed noise  }  \\
\hline
0.1&0.1&   10& 1.57e-3&  2.56e-6& 6.91e-3& 1.30e-3& 2.26e-3& 3.29e-3\\
0.25&0.25& 10& 2.47e-3&  6.01e-6& 1.46e-2& 2.08e-3& 3.62e-3& 5.19e-3\\
0.5 &0.5 &0.1& 1.19e-1&  1.22e-4& 8.74e-1& 7.73e-2& 1.58e-1& 2.784e-1\\
\end{tabular}
\caption{Results for \eqref{Rastrigin} with $d=1$, $N=100,$ $\gamma = 0.1$ and $\xi = 0.0056$}
\label{table:Rastrigin1_new}
\end{table}

We repeated the test for the Rastrigin function with $d=3.$ The parameters of the method and the considered values of $\alpha$ are as in the previous test, while the number of particles is set to $N=500.$ The results are reported in Table \ref{table:Rastrigin3_new}. With these number of particles, in the noisy free case, we manage to approximate the solution with an average error  of the order of $10^{-2}$, and the same accuracy is still obtained in the absolute noise case, while with relative and mixed noise the accuracy is deteriorated for values of $\sigma_1$ grater than 0.1.

\begin{table}[h]
\begin{tabular}{ c |c| c| c| c| c| c| c| c }
$\sigma_0$ & $\sigma_1$ & $\alpha$ & mean err. & min err. & max err. &  50 perc. & 75 perc & 90 perc\\
\hline
\multicolumn{9}{c}{without noise  }  \\
\hline
0 &0   & $1$& 5.07e-02 &  2.04e-03 &  8.57e-01  & 4.375e-02 &  6.19e-02  & 8.18e-02\\
\hline
\multicolumn{9}{c}{absolute noise  }  \\
\hline
0.1 &0 & 1&    4.83e-02&   3.45e-03&   8.25e-01&   4.30e-02  & 58.9e-02  & 8.08e-02
\\
0.25 &0& 1&   4.75e-02&   3.42e-03&   8.50e-01&   4.18e-02&   5.77e-02 &  7.74e-02

\\
0.5 &0 & 1&       5.22e-02 &  5.08e-03&   8.26e-01 &  4.59e-02  & 6.32e-02 &  8.22e-02

\\
1.0 &0 & 1&       5.84e-02 &  5.01e-03 &  8.48e-01&   4.82e-02  & 6.93e-02 &  9.10e-02

\\
\hline
\multicolumn{9}{c}{relative noise  }  \\
\hline
0&0.1&   0.5&       7.43e-02 &   5.20e-03  & 3.58e-01  & 6.86e-02 &  9.52e-02  & 1.21e-01

\\
0&0.25& 0.1&    1.54e-01&   1.72e-02 &  4.31e-01&   1.52e-01   &1.94e-01 &  2.43e-01
\\
0 &0.5 &0.1&    2.48e-01&   9.36e-02 &  8.91e-01 &  2.35e-01 &  3.16e-01 &  4.00e-01\\

\hline
\multicolumn{9}{c}{mixed noise  }  \\
\hline
0.1&0.1&   0.5&       7.44e-02 &   7.11e-03  & 4.93e-01  & 6.83e-02 &  9.36e-02  & 1.22e-01
\\
0.25&0.25& 0.1&    1.50e-01&   9.39e-03 &  4.51e-01&   1.43e-01   &1.87e-01 &  2.40e-01
\\
0.5 &0.5 &0.1&    2.54e-01&   1.31e-02 &  9.79e-01 &  2.37e-01 &  3.22e-01 &  4.15e-01
\end{tabular}
\caption{Results for \eqref{Rastrigin} with $d=3$, $N=500,$ $\gamma = 0.1$ and $\xi = 0.0056$}
\label{table:Rastrigin3_new}
\end{table}

We repeated the test for the rotated Rastrigin function with $d=2.$ That is, we define 
$f(x) = R(Wx)$ where $R$ is the Rastrigin function defined in \eqref{Rastrigin} and $W\in\R^{2\times 2}$ is the rotation matrix with angle $\pi/3.$ In this way we obtain an objective function that, differently from the original Rastrigin function, is not separable and we can then show the effect of noise on non separable functions. 
The number of particles is set to $N=200,$ and the values of $\alpha$ are the same as in the previous tests. The results are reported 
in Table \ref{table:Rastrigin2_rot_new}.  {For the case of absolute noise, the results are analogous to those that we observed for  the Rastrigin function with $d=1$ (Table \ref{table:Rastrigin1_new}), with the method exhibiting significant robustness with respect to noise. On the other hand, in case of relative noise, we notice a higher decrease in the average accuracy even for small values of $\sigma_1.$  }\\

\begin{table}
\begin{tabular}{ c |c| c| c| c| c| c| c| c }
$\sigma_0$ & $\sigma_1$ & $\alpha$ & mean err. & min err. & max err. &  50 perc. & 75 perc & 90 perc\\
\hline
\multicolumn{9}{c}{without noise  }  \\
\hline
0 &0   & $10^4$& 2.18e-04 &  4.37e-06 &  9.73e-04  & 1.95e-04 &  2.68e-042  & 3.91e-04\\
\hline
\multicolumn{9}{c}{absolute noise  }  \\
\hline
0.1 &0 & 1000&    5.95e-03&   9.90e-05&   2.51e-02&   5.23e-03   & 7.15e-03  & 1.07e-02
\\
0.25 &0& 10&   8.93e-03 &  3.16e-04 &  4.08e-02 &  7.69e-03 &  1.05e-02 &  1.65e-02

\\
0.5 &0 & 10&   2.03e-02 &  1.17e-04 &  9.65e-01 &  1.45e-02 &  1.97e-02  & 3.15e-02

\\
1.0 &0 & 1&   3.57e-02 &  1.65e-04 &  8.31e-01  & 3.00e-02 &  4.14e-02 &  6.05e-02

\\
\hline
\multicolumn{9}{c}{relative noise  }  \\
\hline
0&0.1&   1&   4.36e-02 &  1.50e-03  & 9.29e-01  & 3.37e-02  & 4.79e-02  & 7.17e-02

\\
0&0.25& 0.1&   1.65e-01 &  8.29e-03  & 6.10e-01 &  1.44e-01 &  1.97e-01 &  2.98e-01

\\
0 &0.5 &0.1&       2.20e-01  & 5.25e-03 &  7.31e-01 &  2.00e-01  & 2.67e-01   &4.06e-01
\\

\hline
\multicolumn{9}{c}{mixed noise  }  \\
\hline
0.1&0.1&   1&         4.24e-02 &  1.30e-03 &  8.74e-01 &  3.32e-02 &  4.56e-02 &  7.17e-02

\\
0.25&0.25& 0.1&       1.66e-01 &  5.15e-03  & 7.87e-01 &  1.46e-01 &  1.94e-01   &2.95e-01

\\
0.5 &0.5 &0.1&       2.19e-01 &  1.09e-03  & 7.42e-01 &  1.99e-01 &  2.69e-01  & 4.01e-01

\end{tabular}

\captionof{table}{Results for the rotated Rastrigin function with $d=2$, $N=200,$ $\gamma = 0.1$ and $\xi = 0.0056$}
\label{table:Rastrigin2_rot_new}
\end{table}

The condition $\ca<1$ can be weakened employing node-independent random diffusion terms, i.e, $\eta_{k,s}^i = \eta_{k,s}$ for every $i=1,\dots,N$ as in 
\cite{discreteCBO} where convergence for the noise-free case is proved 
assuming $(1-\gamma)^2+\xi^2<1$.  On the one hand, employing different  diffusion terms for each particle enhance the exploration capability of the particles.  On the other hand, the condition $(1-\gamma)^2+\xi^2<1$ is satisfied for larger values of $\xi$, which is also known to increase the exploratory power of the system of particles.
To give more insight into the method's characteristic, we repeated the  tests above, with parameters $\gamma=0.01$ and $\xi = 0.1.$ We remark that while this choice of parameters does not guarantee convergence according to our analysis, convergence was achieved in practice in all our runs, suggesting that the convergence condition may be too conservative. In  Table \ref{table:Rastrigin3old} we report the results obtained  for $d=3$ and $N=100$ particles, showing that this relaxed choice of the parameters allows to  use  a smaller number of particles. 
\begin{table}
\begin{tabular}{ c |c| c| c| c| c| c| c| c }
$\sigma_0$ & $\sigma_1$ & $\alpha$ & mean err. & min err. & max err. &  50 perc. & 75 perc & 90 perc\\
\hline
\multicolumn{9}{c}{without noise  }  \\
\hline
0 &0   & $10^4$& 1.07e-02 &  1.34e-04 &  9.95e-01  & 1.41e-03 &  2.20e-03  & 4.40e-03\\
\hline
\multicolumn{9}{c}{absolute noise  }  \\
\hline
0.1 &0 & 1000&    3.51e-02&   1.58e-03&   9.95e-01&   1.35e-02  & 1.92e-02  & 2.62e-02
\\
0.25 &0& 1000&   8.82e-02&   3.94e-03&   1.41e+00&   1.76e-02&   2.52e-02 &  3.88e-02

\\
0.5 &0 & 1000&       1.18e-01 &  1.64e-03&   1.40e+00 &  2.15e-02  & 3.06e-02 &  5.59e-02

\\
1.0 &0 & 100&       9.24e-02 &  2.05e-03 &  1.40e+00&   2.60e-02  & 3.56e-02 &  5.01e-02

\\

\hline
\multicolumn{9}{c}{relative noise  }  \\
\hline
0&0.1&   1000&       1.05e-01 &   5.31e-03  & 1.36e+00  & 4.10e-02 &  5.65e-02  & 8.37e-02

\\
0&0.25& 0.1&    4.41e-01&   5.89e-03 &  1.35e+00&   4.13e-02   &5.74e-01 &  7.46e-02
\\
0 &0.5 &0.1&    8.30e-01&   7.18e-02 &  7.79e+01 &  6.13e-01 &  8.73e-01 &  1.14e00\\

\hline
\multicolumn{9}{c}{mixed noise  }  \\
\hline
0.1&0.1&   1000&       1.01e-01 &   2.44e-03  & 1.13e+00  & 4.01e-02 &  5.50e-02  & 7.87e-02
\\
0.25&0.25& 0.1&    4.50e-01&   3.80e-02 &  1.27e+00&   4.28e-01   &5.82e-01 &  7.40e-01
\\
0.5 &0.5 &0.1&    6.92e-01&   4.16e-02 &  5.75e+00 &  6.17e-01 &  8.36e-01 &  1.11e+00

\end{tabular}

\captionof{table}{Results for \eqref{Rastrigin} with $d=3$, $N=100,$ $\gamma = 0.01$ and $\xi = 0.1$.}
\label{table:Rastrigin3old}
\end{table}

\subsection{Subsampling}

We now consider a finite-sum minimization problem and we combine the CBO method with a subsampling strategy. That is, we consider the case where the objective function is of the form 
\be \label{finite-sum}
f(\xbf) = \frac{1}{M}\sum_{j=1}^M f_j(\xbf)\enspace \text{with}\ f_j:\R^d\longrightarrow\R
\ee
and the oracle $\hat f$ is define as
\begin{equation}\label{oracle_sampling}
\hat f(\xbf) = \frac{1}{|S|}\sum_{j\in S} f_j(\xbf)\end{equation}
where $S$ is a random subset of $\{1,\dots,M\}$ chosen randomly and uniformly. We recall that $\E_k[\hat f(\xbf_k^i)]=f(\xbf_k^i)$.
We consider the sigmoid function coupled with the quadratic loss  for a binary classification problem. That is, given a dataset $\{(a_j, b_j)\}_{j=1}^M$, with $a_j\in\R^d$ and $b_j\in\{0,1\}$,
$f_j$ in \eqref{finite-sum} takes the form
 \begin{equation}\label{class_objfun} f_j(\xbf) = \lr{b_j - \frac{1}{1+e^{-\xbf^\top a_j }}}^2.\end{equation}
We apply the CBO method given in \eqref{CBOhat}, with $\hat f(\xbf^i_k)$ as in \eqref{oracle_sampling}.  Specifically, at very iteration $k$ and for every particle $i$ we generate $S^i_k$ random subset of $\{1,\dots,M\}$ with uniform probability and cardinality $|S^i_k| = \ell M$, for a given $\ell \in(0,1),$ and we define
$$\hat f(\xbf^i_k) = \frac{1}{|S^i_k|}\sum_{j\in S^i_k} f_j(\xbf^i_k).$$
We consider different values of the sampling parameter $\ell$ and compare the performance of the resulting methods, in terms of accuracy and computational cost. In the following, we report the results of the experiments on the \emph{Rice (Cammeo and Osmancik)} dataset \cite{rice}, which has $7$ features and $3810$ instances. We split the dataset into a training set and a testing set, with $2857$ and $953$ instances respectively. The training set is used in the definition of the objective function \eqref{class_objfun}, so we have $d = 7$ and $M = 2857$. For each value of $\ell=1,\ 0.75,\ 0.5,\ 0.25\ 0.1,\ 0.05,\ 0.025,\ 0.01$, we run the CBO method 100 times, with $N=1000$ particles and parameters $\gamma = 0.1$, $\xi = 0.0056$, and $\alpha = 10^3.$ For each instance of the method, the initial position $\xbf^i_0$ of every particle is chosen randomly on the set $[-10^3,\ 10^3]^d$ with uniform distribution, while the execution is terminate when 
$$\frac{1}{N}\sum_{i=1}^N\|\xbf^i_k-\xbf^{av}_k\|\leq 10^{-3},$$
where $\xbf^{av}_k$ denotes the average of the vectors $\xbf^i_k$ for $i=1,\dots,N$.
To evaluate the performance of the method we consider the accuracy and the computational cost. The accuracy is computed as the percentage of instances in the testing set that are correctly classified, using as vector of parameters the consensus point $\xbf^\alpha_k$ at termination.  The results are summarized in Table \ref{table:classification_new}
where for every  value of the sampling parameter $\ell$ we report the average accuracy over all runs (\emph{mean acc.}), the average number of iterations (\emph{mean it}), the average number of function $f_j$ evaluations divided by $M$ (\emph{mean eval.}) and the average computational cost (\emph{mean cost}). The computational cost is computed as $\bar k\,d(\lceil \ell\,M\rceil +2),$ where $\bar k$ is the number of iterations and we assigned cost $d\,\lceil \ell\,M\rceil$ to the evaluation of $\hat f$ and cost $2d$ to the update of a particle's position.

\begin{figure}[h]
    \centering
    \includegraphics[width=0.60\textwidth]{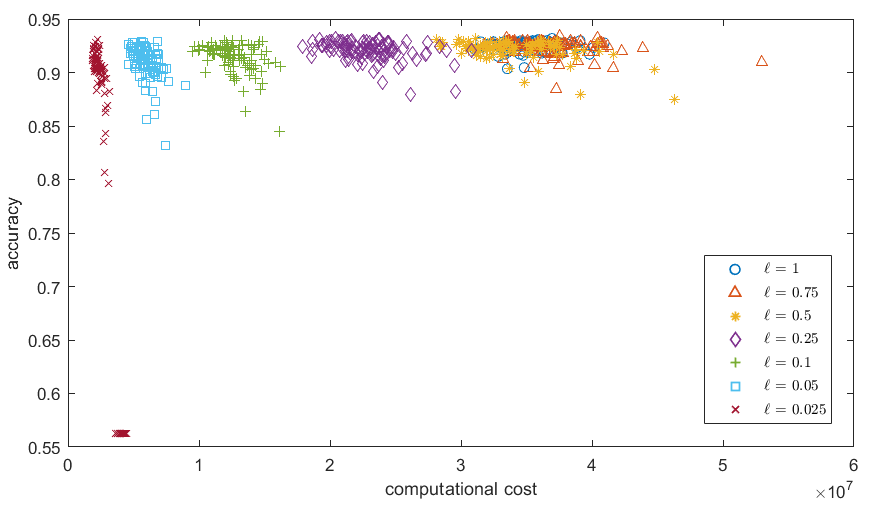}
    \caption{Accuracy vs computational cost on problem \eqref{class_objfun} on the Rice dataset for different values of $\ell$, for $N=500,$ $\gamma = 0.01$ and $\xi = 0.1$}
    \label{fig:classification}
\end{figure}

We repeated the tests with parameters $\gamma = 0.01$ and $\xi = 0.1.$ In this case, we reduced the number of particles to $N=500.$ In table \ref{table:classification} we report the obtained results, while in figure \ref{fig:classification} we show the
accuracy versus the 
computational cost; we plot a point for each of the 100 runs that we executed for every value of the sampling parameter $\ell$.

\begin{table}[h]
\centering
\begin{tabular}{ c| c| c| c| c }
 $\ell$ & mean acc. & mean it & mean eval. &mean cost\\
 \hline
 1 & 82.5\% & 136.7 & 1.36730e+5 & 2.7364e+6\\
 0.75 & 83.4\% & 136.7 & 1.0250e+5 & 2.0518e+6\\
 0.5 & 82.1\% & 136.7 & 6.8374e+4 & 1.3693e+6\\
 0.25 & 81.3\% & 136.7 & 3.4221e+4 & 6.8630e+5\\
 0.1 & 78.5\% & 136.7 & 1.3692e+4 & 2.7575e+5\\
 0.05 & 75.1\% & 136.7 & 6.8417e+3 & 1.3874e+5\\
 0.025 & 66.1\% & 136.7 & 3.4453e+3 & 7.0816e+4\\
\end{tabular}
\caption{Summary of the experiments on problem \eqref{class_objfun} on the Rice dataset for different values of $\ell$, with $N=1000$, $\gamma=0.1$ and $\xi = 0.0056$}
\label{table:classification_new}
\end{table}

\begin{table}[h]
\centering
\begin{tabular}{ c| c| c| c| c }
 $\ell$ & mean acc. & mean it & mean eval. &mean cost\\
 \hline
 1 & 92.5\% & 1745.3 & 8.7263e+5 & 3.4928e+7\\
 0.75 & 92.3\% & 2456.9 & 9.2135e+5 & 3.6891e+7\\
 0.5 & 92.1\% & 3392.9 & 8.4823e+5 & 3.3987e+7\\
 0.25 & 92.0\% & 4578.8 & 5.7234e+5 & 2.2981e+7\\
 0.1 & 91.4\% & 6283.5 & 3.1417e+5 & 1.2668e+7\\
 0.05 & 91.0\% & 5828.1 & 1.4563e+5 & 5.9125e+6\\
 0.025 & 85.0\% & 5003.5 & 6.2544e+4 & 2.5918e+6\\
\end{tabular}
\caption{Summary of the experiments on problem \eqref{class_objfun} on the Rice dataset for different values of $\ell$, with $N=500,$ $\gamma=0.01$ and $\xi = 0.1$}
\label{table:classification}
\end{table}

For both choices of the parameters we can observe that the final accuracy is not greatly affected by the subsampling, up to $\ell=0.25$ in the first case and $\ell=0.05$ in the second case. 
For $\gamma=0.1$ and $\xi = 0.0056$ (Table \ref{table:classification_new}) the number of iterations is approximately constant for all considered values of $\ell$, which implies that the computational cost decreases significantly as the sampling size decreases. For $\ell=0.5$ the overall computational cost is halved, but the average accuracy is less than $0.5\%$ smaller.
For $\gamma=0.01$ and $\xi = 0.1$ (Table \ref{table:classification}) the number of iterations increases as the sampling factor decreases, but the overall computational cost decreases. The method that employs $\ell=0.1$ is approximately 2.5 times less expensive than the method that uses the full sample, while the average accuracy is only around $1\%$ smaller.
Comparing the results for the two choices of the parameters (namely $\gamma = 0.1,\ \xi = 0.0056$ and $\gamma = 0.01,\ \xi = 0.1$) we notice that, despite the fact that the number of particles $N$ is larger in the first case, the particles converge in a smaller number of iterations, and that the overall accuracy is lower, compared to that reached by the method that employs smaller $\gamma$ and larger $\xi.$ This confirms that the condition $\ca<1$ from the theoretical analysis may be too restrictive and the method can perform better in practice for choices of the parameters that increase the exploration of the search space.

Figure \ref{fig:classification} clearly depicts the gain in terms of computational cost with $\ell\le 0.25$ and suggests that $\ell=0.25$ is a good compromise between accuracy and computational cost.

{\section{Conclusions}\label{sec:conclusions}
We considered a discrete-time CBO method with particle-independent diffusion noises in the case where only a stochastic estimator of the objective function is available for evaluation. We proved that, for suitable choices of the parameters of the methods, the particles still achieve asymptotic consensus, and we study the optimality gap at the consensus point. Moreover, given any accuracy threshold, we provide conditions on the objective function and the noise of the stochastic oracle that ensure that the mean-squared error becomes, in expectation, smaller than the chosen threshold. 
We studied numerically the influence of noise on the performance of the method and showed that, when the exact objective function is not available, the method can still achieve good results.}

\section*{Acknowledgements}
 The research that led to the present paper was partially supported by INDAM-GNCS through Progetti di Ricerca 2024 and 
by PNRR - Missione 4 Istruzione e Ricerca - Componente C2 Investimento 1.1, Fondo per il Programma Nazionale di Ricerca e Progetti di Rilevante Interesse Nazionale (PRIN) funded by the European Commission under the NextGeneration EU programme, project ``Advanced optimization METhods for automated central veIn Sign detection in multiple sclerosis from magneTic resonAnce imaging (AMETISTA)'',  code: P2022J9SNP,
MUR D.D. financing decree n. 1379 of 1st September 2023 (CUP E53D23017980001), project 
``Numerical Optimization with Adaptive Accuracy and Applications to Machine Learning'',  code: 2022N3ZNAX
 MUR D.D. financing decree n. 973 of 30th June 2023 (CUP B53D23012670006).

\appendix 
\section*{Appendix A}

    \begin{lemma}\label{lemma:gbounds}\cite{gaussianbound}
        Let $\{v_1,\dots v_n\}_{i=1}^N$ be Gaussian random variables with $v_i\sim\mathcal{N}(0,\xi_i^2)$ for every $i.$ Then
        $$\E\lrq{\max_{i=1:N}v_i^2}\leq 4\log(\sqrt{2}N)\max_{i=1:N}\xi_i^2.$$
    \end{lemma}
For the proof of Lemma \ref{lemma_ergM} we will use the following results.

\begin{lemma}\cite[Lemma 3.3 ]{hetCBO}\label{lemma_erg}
If $M\in\R^{N\times N}$ is such that $\sum_{j=1}^N M_{ij} = m$ for any index $i=1,\dots,N$, then for every vector $\ybf\in\R^N$ we have 
$\diam(My)\leq (m-\erg(M))\diam(y).$
\end{lemma}
\begin{lemma}\cite[Lemma 4.1]{hetCBO}\label{lemma_erg2}
Given any matrix $M\in\R^{N\times N}$, the ergodicity coefficient satisfies the following inequality: 
$\erg(M)\geq -2\|M\|_{\infty}$ for any $M\in\R^{N\times N}$
\end{lemma}

\begin{proof} [Proof of Lemma \ref{lemma_ergM}]
First, notice that for any $i,j=1,\dots,N$ we have 
\be \label{Mks}
(M^s_k)_{ij} = \begin{cases}1-\gamma -  \eta^i_{k,s} + (\gamma + \eta^i_{k,s})\hat v^\alpha_{k,j} & \text{if } i=j\\
(\gamma + \eta^i_{k,s})\hat v^\alpha_{k,j} & \text{if } i\ne j.
\end{cases}
\ee
Therefore, for any $i=1,\dots, N$ we have 
$$\sum_{j=1}^N M_{ij} = 1-\gamma -  \eta^i_{k,s} + (\gamma + \eta^i_{k,s})\hat v^\alpha_{k,i} + (\gamma + \eta^i_{k,s})\sum_{j\neq i} \hat v^\alpha_{k,j} = 1,$$
where we used the fact that, by definition,
$$\sum_{j=1}^N \hat v^\alpha_{k,j} = \sum_{j=1}^N\frac{e^{-\alpha\hat f(\xbf^j_k, \omega^j_k)}}{\sum_{h=1}^N e^{-\alpha\hat f(\xbf^h_k, \omega^h_k)}} = 1.$$
Then, the thesis in $i)$ directly follows by  \eqref{Y2} and Lemma \ref{lemma_erg}.

Regarding $ii)$,  $M^s_k = A^s_k+B^s_k$ with 
$$A^s_k = (1-\gamma)I_N + \gamma \hat V^\alpha_k\enspace\text{and}\enspace B^s_k =  - D^s_k (I_N-\hat V^\alpha_k), $$
therefore, by 
the definition of ergodicity given in \eqref{ergo} we have
\begin{equation}\label{ergM}\erg(M^s_k)\geq \erg(A^s_k)+ \erg(B^s_k).\end{equation}
We now find a lower bound to the ergodicity coefficient of $A^s_k$ and $B^s_k$.
By definition of $A^s_k$, we have
$$(A^s_k)_{ij} = 1-\gamma + \gamma\hat v^\alpha_{k,i} \quad \text{if } i=j, \quad 
 (A^s_k)_{ij} =\gamma \hat v^\alpha_{k,j}  \quad\text{otherwise}. $$
Since $1-\gamma>0$, this yields
\begin{equation}\label{erg_A}
\begin{aligned}
\erg(A^s_k) &= \min_{i,l=1:N}\sum_{j=1}^N \min\{(A^s_k)_{ij}, (A^s_k)_{lj}\} = 
 \min_{i,l} \sum_{j=1}^N \gamma \hat v^\alpha_{k,j} = \gamma.
\end{aligned}
\end{equation}
Let us now consider $B^s_k.$ By Lemma \ref{lemma_erg2}, we have
\begin{equation*}
\erg(B^s_k)\geq -2\|B^s_k\|_\infty
\geq-2\|D^s_k\|_\infty(\|I_N\|_\infty+\|\hat V^\alpha_k\|_\infty)\geq -2\max_{i=1:N}|\eta^i_{k,s}|(1+\|\hat V^\alpha_k\|_\infty).
\end{equation*}
Since $\sum_{j=1}^N \hat v^\alpha_{kj} = 1$ and $\hat v^\alpha_{kj}>0$ for every $j$, we have that $\|\hat V^\alpha_k\|_\infty = 1$ and thus 
\begin{equation}\label{erg_B}
\erg(B^s_k)\geq -4\max_{i=1:N}|\eta^i_{k,s}|.
\end{equation}

Replacing \eqref{erg_A} and \eqref{erg_B} into \eqref{ergM} we get the thesis.
\qed\end{proof}

\begin{proof}[Proof of Theorem \ref{convergencethm}]
    For every $k\in\N$ let us define 
    $$Q_{k,s} = \gamma\sum_{l=0}^k(x^i_{l,s}-\hat x^\alpha_{l,s}),\enspace P_{k,s} = \sum_{l=0}^k(x^i_{l,s}-\hat x^\alpha_{l,s})\eta_{l,s}^i.$$
    Recursively applying \eqref{CBOhat} we get
    \be\label{AkBk}
    \begin{aligned}
    x^i_{k+1,s} = \xbf^i_{0,s} - \gamma Q_{k,s}-P_{k,s}.
    \end{aligned}
    \ee
    We prove that $Q_{k,s}$ and $P_{k,s}$ converge almost surely as $k$ tends to $+\infty$. Let us first consider $Q_{k,s}.$
    Using the fact that $\sum_{j=1:N}\hat v^\alpha_{kj} = 1$, \eqref{th2_1} and part $ii)$ of Lemma \ref{lemma_diam_CBO}, we have 
$$
    |x^i_{l,s}-\hat x^\alpha_{l,s}| 
    \leq e^{-lS_{l,s}}\diam\lr{\ybf^s_0}.
    $$
    Since $\lim_{i \rightarrow \infty}S_{l,s}\ge 1-\ca$  almost surely,
   there exist two stochastic variables $C_1,C_2>0$ such that
    \be \label{conv_component}
    |x^i_{l,s}-\hat x^\alpha_{l,s}|\leq C_1 e^{-lC_2}\enspace \text{almost surely.}
    \ee
    Consider now the sequence $\{\widetilde Q_{k,s}\}$ defined as
    $$\widetilde Q_{k,s} = Q_{k,s} - \gamma \sum_{l=0}^{k}C_1e^{-lC_2} = \gamma \sum_{l=0}^{k}\lr{x^i_{l,s}-\hat x^\alpha_{l,s} - C_1e^{-lC_2}}. $$
    The sequence $\{\tilde Q_{k,s}\}$ is non-increasing by \eqref{conv_component}. Moreover
    \be
    \begin{aligned}
    \widetilde Q_{k,s} &= Q_{k,s} + \gamma\sum_{l=0}^{k}C_1e^{-lC_2} - 2\gamma\sum_{l=0}^{k}C_1e^{-lC_2}    \\  
    &=\gamma\sum_{l=0}^{k}\lr{x^i_{l,s}-\hat x^\alpha_{l,s} + C_1e^{-lC_2}} - 2\gamma\sum_{l=0}^{k}C_1e^{-lC_2}  \\  
    &\geq - 2\gamma C_1\sum_{l=0}^{k}e^{-lC_2} = -2\gamma C_1 \frac{1-e^{C_2(k+1)}}{1-e^{C_2}}\geq \frac{1}{1-e^{C_2}}.
    \end{aligned}
    \ee
    Since almost surely we have that $\{\widetilde Q_{k,s}\}$ is a non-increasing sequence, bounded from below, we have that there exists $\widetilde Q_{\infty,s}\in\R$ such that
    $\lim_{k\rightarrow+\infty}\widetilde Q_{k,s} = \widetilde Q_{\infty,s}\enspace \text{almost surely},$
    which in turn implies that, almost surely,
    $$\lim_{k\rightarrow+\infty} Q_{k,s} = \lim_{k\rightarrow+\infty} \lr{\widetilde Q_{k,s} + \gamma \sum_{l=0}^{k}C_1e^{-lC_2} }= \widetilde Q_{\infty,s} + \gamma C_1 \frac{1}{1-e^{-C_2}} =:Q_{\infty,s}.$$
    Let us now consider the sequence $\{P_{k,s}\}.$
    By the independence of $\eta_{k,s}^i$ and $(x^i_{k,s}-\hat x^\alpha_{k,s})$, and the fact that $\E[\eta_{k,s}^i] = 0$, we have
$$
    \begin{aligned}
    & E[P_{k,s}|P_{i,s}\ i=0:k-1] = \E\lrq{\sum_{l=0}^k(x^i_{l,s}-\hat x^\alpha_{l,s})\eta_{l,s}^i |P_{i,s}\ i=0:k-1}   \\  
    &= \E\lrq{P_{k-1,s} +(x^i_{k,s}-\hat x^\alpha_{k,s})\eta_{k,s}^i |P_{i,s}\ i=0:k-1}   \\  
    & = P_{k-1,s} + \E\lrq{(x^i_{k,s}-\hat x^\alpha_{k,s})\eta_{k,s}^i |P_{i,s}\ i=0:k-1}   \\  
    & = P_{k-1,s} + \E[\eta_{k,s}^i]\E\lrq{(x^i_{k,s}-\hat x^\alpha_{k,s})|P_{i,s}\ i=0:k-1} = P_{k-1,s}.
    \end{aligned}
    $$
    That is, $P_{k,s}$ is a Martingale.
    Moreover we have
    $$
    \begin{aligned}
    &\E[P_{k,s}^2] = \E\lrq{\lr{\sum_{l=0}^k(x^i_{l,s}-\hat x^\alpha_{l,s})\eta_{l,s}^i}^2}  \\  & = \E\lrq{\sum_{l=0}^k(x^i_{l,s}-\hat x^\alpha_{l,s})^2(\eta_{l,s}^i)^2 + \sum_{l=0}^{k-1}\sum_{\substack{p=0 \\ p\neq l}}^{k-1}(x^i_{l,s}-\hat x^\alpha_{l,s})(x^i_{p,s}-\hat x^\alpha_{p,s})\eta_{l,s}^i\eta_{p,s}^i }.
     \end{aligned}
    $$
     Letting $q_{max}=\max\{l,p\}$ and $q_{min}=\min\{l,p\}$ we note
  that $\eta_{q_{max},s}$ is independent of $(x^i_{l,s}-\hat x^\alpha_{l,s})(x^i_{p,s}-\hat x^\alpha_{p,s})\eta_{q_{min},s}^i.
  $
  Then, by \eqref{th2_1} and Lemma \ref{lemma_diam_CBO} we get
  \be\label{expB2}
    \begin{aligned}
    \E[P_{k,s}^2] &=
\xi^2\E\lrq{\sum_{l=0}^k(x^i_{l,s}-\hat x^\alpha_{l,s})^2}  \\   &+\sum_{l=0}^{k-1}\sum_{\substack{p=0 \\ p\neq l}}^{k-1}\E[(x^i_{l,s}-\hat x^\alpha_{l,s})(x^i_{p,s}-\hat x^\alpha_{p,s})\eta_{q_{min},s}^i]\E[\eta_{q_{max},s}^i]   \\  
    & = \xi^2\E\lrq{\sum_{l=0}^k(x^i_{l,s}-\hat x^\alpha_{l,s})^2}\leq \xi^2\sum_{l=0}^k\E\lrq{\diam(\ybf^s_l)^2}\\
    &\leq \xi^2\E\lrq{\diam(\ybf^s_0)^2}\sum_{l=0}^k\E\lrq{e^{-2lS_{l,s}}}.
    \end{aligned}
    \ee
    Since $\lim_{l\rightarrow+\infty}S_{l,s} \ge  1-\ca>0$ almost surely, then there exists a constant $\bar{S}>0$ such that for $l$ large enough $S_{l,s}>\bar{S}$ and thus $e^{-lS_{l,s}}\leq e^{-l\bar{S}}.$
    This, together with \eqref{expB2} ensures that $\E[P_{k,s}^2]$ is uniformly bounded and therefore, by Doob's theorem, there exists $P_{\infty,s}$ such that 
    $\lim_{k\rightarrow+\infty}P_{k,s} = P_{\infty,s}\enspace \text{almost surely}.$
    Finally, from \eqref{AkBk} and the convergence of $Q_{k,s}$ and $P_{k,s}$ we have, almost surely  for every $s = 1,\dots,d$,
    \be
    \begin{aligned}
\lim_{k\rightarrow+\infty}x^i_{k+1,s} = \lim_{k\rightarrow+\infty}(\xbf^i_{0,s} - \gamma Q_{k,s}-P_{k,s}) = \xbf^i_{0,s} - \gamma Q_{\infty,s}-P_{\infty,s}.
    \end{aligned}
    \ee
   Then,  for every index $i$ there exists $\xbf^i_{\infty}$ such that 
$$\lim_{k\rightarrow+\infty}\xbf^i_k = \xbf^i_{\infty}\enspace \text{almost surely}.$$
    Since Th. \ref{thm_cons1} yields
$\lim_{k\rightarrow+\infty}\|\xbf^i_k - \xbf^j_k \|= 0$ a. s.  for every $i,j=1,\dots,N$,
    we conclude that there exists $\xbf_\infty\in\R^d$ such that 
    $\xbf^i_\infty = \xbf_\infty$ for every $i=1,\dots,N$, and  the thesis follows.
\qed\end{proof}

\begin{proof}[Proof of Lemma \ref{lemmaAk}]
By \eqref{CBOhat_compatto} it follows 
 \be\label{consecutivediff}
x^i_{k+1}-x^i_{k} = - (\gamma I_d+\diageta_{k}^i)\lr{x^i_{k}-\hat x^\alpha_{k}},
\ee
    and similarly
    \be\label{thetadiff}(\theta^i_k\xbf^i_{k+1}+(1-\theta^i_k)\xbf^i_{k}-\hat\xbf^\alpha_{k})= 
    ((1-\theta^i_k\gamma) I_d-\theta^i_k\diageta_{k}^i)\lr{x^i_{k}-\hat x^\alpha_{k}}.
    \ee
    Let  $\psi(\theta) = (1-\theta\gamma-\theta\eta_{k,s}^i)^2$, 
    since $\psi(\theta)$ is a convex quadratic function, we have that
   \be\label{psibound}\max_{\theta\in(0,1)}\psi(\theta) = \max\{\psi(0),\psi(1)\}\leq \psi(0)+\psi(1) = 1+(1-\gamma-\eta_{k,s}^i)^2.\ee
   Then, 
    using the Cauchy-Schwartz inequality, the Lipschitz continuity of the gradient, \eqref{consecutivediff},\eqref{thetadiff} and \eqref{psibound} we  have
    \be
    \begin{aligned}
    A_k 
    &\leq \frac{M_H}{N}e^{-\alpha f_*}\sum_{i=1}^N\lr{\sum_{s=1}^d(1+(1-\gamma-\eta_{k,s}^i)^2(x^i_{k,s}-\hat x^\alpha_{k,s})^2}^{1/2}\cdot\\
&\enspace\cdot\lr{\sum_{s=1}^d(\gamma+\eta_{k,s}^i)^2(x^i_{k,s}-\hat x^\alpha_{k,s})^2}^{1/2}.
    \end{aligned}
    \ee
    Taking the expected value on both sides of the inequality above, using Holder inequality, the linearity of $\E[\ \cdot\ ]$ and the independence of $\eta_{k,s}$ and $(x^i_{k,s}-\hat x^\alpha_{k,s})$ we have
   \bes
    \begin{aligned}
    \E[A_k]  
    &\leq \frac{M_H}{N}e^{-\alpha f_*}\sum_{i=1}^N\bigg(\sum_{s=1}^d(1+(1-\gamma)^2+\xi^2)\E[(x^i_{k,s}-\hat x^\alpha_{k,s})^2]\bigg)^{1/2}\cdot\\
    &\enspace\cdot \bigg(\sum_{s=1}^d(\gamma^2+\xi^2)\E[(x^i_{k,s}-\hat x^\alpha_{k,s})^2]\bigg)^{1/2}   \\  
    &=\frac{M_H}{N}(1+(1-\gamma)^2+\xi^2)^{1/2}(\gamma^2+\xi^2)^{1/2}e^{-\alpha f_*}\sum_{i=1}^N\E[\|\xbf^i_k-\hat\xbf^\alpha_k\|^2].
    \end{aligned}
    \ees
    By definition of $\Gamma_A$ and applying  (\emph{iv}) in Lemma \ref{lemma_diam_CBO}, we get the thesis.
    \qed\end{proof}

\begin{proof}[Proof of Lemma \ref{diffxalphasol}] 
We follow the proof of Proposition 21 in \cite{Fornasier1}. 
Let $w_{i,k}^\alpha =  e^{-\alpha f(\xbf_k^i)}.$ 
By definition of $\xbf^\alpha_k$ given in \eqref{xkalpha} 
and the triangle inequality, for any $\tilde r\geq r$ we have
\be\label{diffsol}\begin{aligned}
\|\xbf^\alpha_k-\xbf^*\|& = \left\|\frac{1}{\sum_{i=1}^N w_{i,k}^\alpha}\sum_{i=1}^N w^\alpha_{i,k}(\xbf_k^i-\xbf^*)  \right\|  
\le \frac{1}{\sum_{i=1}^N w_{i,k}^\alpha}\sum_{i=1}^N w^\alpha_{i,k}\left \|\xbf_k^i-\xbf^*\right \|    \\ 
&=\frac{1}{\sum_{i=1}^N w_{i,k}^\alpha}\left ( \sum_{i \in {\cal I}_{k,\tilde r}} w^\alpha_{i,k}  \left \|\xbf_k^i-\xbf^*\right \|  
+ \sum_{i \not\in {\cal I}_{k,\tilde r}} w^\alpha_{i,k}   \left \|\xbf_k^i-\xbf^*\right \|   \right ) \\  
&\leq \tilde r    + \frac{\max_{i \not\in {\cal I}_{k,\tilde r}} w_{i,k}^\alpha }{\sum_{i=1}^N w_{i,k}^\alpha}
 \sum_{i \not\in {\cal I}_{k,\tilde r}}  \left \|\xbf_k^i-\xbf^*\right \| \\  
&\leq\tilde r + \frac{\exp(-\alpha\min_{\xbf_k^i \in \B_{\tilde r}^\complement}f(\xbf_k^i))}
{\sum_{i=1}^N w_{i,k}^\alpha}  \sum_{i \not\in {\cal I}_{k,\tilde r}} w^\alpha_{i,k}   \left \|\xbf_k^i-\xbf^*\right \|.
\end{aligned}
\ee

We take $\tilde r=\frac{1}{\beta}(q+f_r)^{\nu}.$
Notice that $\tilde r\ge r$  as by the assumption on $q$ and the fact that $f_*\geq 0,$ we have
$$\tilde r = \frac{1}{\beta}(q+f_r)^{\nu} \ge \frac{1}{\beta}(f_r-f_*)^\nu=\max_{x\in \B_r}\frac{1}{\beta}(f(x)-f_*)^\nu
\geq \max_{x\in \B_r}\|\xbf-\xbf^*\| = r.$$
By assumption A6 and the definition of $\tilde r$ we have 
\be\label{inff} \min_{\xbf_k^i \in \B_{\tilde r}^\complement}f(\xbf_k^i))\geq\min\left\{f_*+f_{\infty}, f_*+(\beta\tilde r)^{1/\nu}\right\} \geq (\beta\tilde r)^{1/\nu} = q+f_r,\ee
 where we used the fact that $(\beta{\tilde r})^{1/\nu} = q+f_r\leq f_{\infty}.$
 We now distinguish between case $a)$ and $b)$ to bound $\sum_{i=1}^N w_{i,k}^\alpha$.
Let us consider case $a).$
We have 
$$
\sum_{i=1}^N w_{i,k}^\alpha=\sum_{i=1}^N  e^{-\alpha f(\xbf_k^i)}\ge \sum_{\xbf_k^i \in \B_{ r}}  e^{-\alpha f(\xbf_k^i)}\ge |{\cal I}_{k,r}| e^{-\alpha f_r}.
$$
Using this inequality, \eqref{inff} and the definition of $\tilde r$ into \eqref{diffsol} we get the thesis.\\
In case $b)$ we can proceed as above noting that the following inequality holds:
$$
\sum_{i=1}^N w_{i,k}^\alpha=\sum_{i=1}^N  e^{-\alpha f(\xbf_k^i)}\ge
N  e^{-\alpha \max_{i=1:N} f(\xbf^i_k)}.$$
\qed\end{proof}

\begin{proof}[Proof of Theorem \ref{thmVkprob}] 
We follow the proof of Theorem \ref{thmVkdet}. In particular,  by
\eqref{Vkmu}  we have that 
\be\E[V_{k+1}]\leq\mu\E[V_k],\ee
for every $k=0,\dots,\bar k-1$ where $\bar k $ is such that \eqref{bark} holds. 
This  implies the thesis if $\bar k>k_*$ or if $\bar k\le k_*$ and $\E[V_k]\leq\varepsilon$. Therefore the only thing left to show is that the remaining case, $\bar k\le k_*$ and $\E[\|\xbf^\alpha_{\bar k}-\xbf^*\|^2]>C^2_{\bar k}$, cannot happen.
Proceeding as in \eqref{eq:deriv} and using Lemma \ref{diffxalphasol} we have
\be\label{newsplit}\begin{aligned}
    &\|\hat x^\alpha_{\bar k}-x^*\|^2\leq \|x^\alpha_{\bar k}-x^*\|^2+\|\hat x^\alpha_{\bar k}-x^\alpha_{\bar k}\|^2+2\|x^\alpha_{\bar k}-x^*\|\|\hat x^\alpha_{\bar k}-x^\alpha_{\bar k}\| \\
    &\leq \frac{1}{4}C_{\bar k} + N^2V_{\bar k}E_{\bar k}(\alpha)^2 + C_{\bar k} N V_{\bar k}^{1/2}E_{\bar k}(\alpha) \\
    & + \alpha^2\|\mathcal{E}_k\|_{\infty}^2\lr{\sum_{s=1}^d\diam(\ybf^s_0)^2\prod_{l=0}^{{\bar k}}(1-\gamma+4\max_{i=1:N}|\eta_{l,s}^i|)^2} \\
    & + C_{\bar k}\alpha\|\mathcal{E}_k\|_{\infty}\lr{\sum_{s=1}^d\diam(\ybf^s_0)^2\prod_{l=0}^{{\bar k}}(1-\gamma+4\max_{i=1:N}|\eta_{l,s}^i|)^2}^{1/2} \\
    &+ 2\alpha N\|\mathcal{E}_k\|_{\infty}\lr{\sum_{s=1}^d\diam(\ybf^s_0)^2\prod_{l=0}^{{\bar k}}(1-\gamma+4\max_{i=1:N}|\eta_{l,s}^i|)^2}^{1/2}V_{\bar k}^{1/2} E_{\bar k}(\alpha),
\end{aligned}
\ee
where $E_{\bar k}(\alpha)$ comes from Lemma \ref{diffxalphasol} and is defined as
$$E_{\bar k}(\alpha) = \begin{cases}
    e^{-\alpha q} & \text{if } |I_{\bar k, r}|>0\\
    \frac{e^{-\alpha q}}{Ne^{-\alpha\max_{i=1:N} f(\xbf^i_k)}} & \text{otherwise}.
\end{cases}$$
Taking the expected value in \eqref{newsplit} and using Holder's inequality we get
\be \label{newsplit2}\begin{aligned}
    &\E\lrq{\|\hat x^\alpha_{\bar k}-x^*\|^2}\leq \frac{1}{4}C_{\bar k} + N^2\E\lrq{V_{\bar k}^2}^{1/2}\E\lrq{E_{\bar k}(\alpha)^4}^{1/2} + C_{\bar k} N \E\lrq{V_{\bar k}}^{1/2}\E\lrq{E_{\bar k}(\alpha)^2}^{1/2} \\
    & + \alpha^2M_v^2\ca^{2\bar k}D_0 
     + C_{\bar k}\alpha M_v \ca^{\bar k}D_0^{1/2} 
    + 2\alpha N M_v \ca^{\bar k}D_0^{1/2}\E\lrq{V_{\bar k}^2}^{1/4} \E\lrq{E_{\bar k}(\alpha)^4}^{1/4}.
\end{aligned}
\ee
By the law of total expectation we have that, for any $\ell\in\N$
$$
    \begin{aligned}
        \E\lrq{E_{\bar k}(\alpha)^\ell} = &\E[E_{\bar k}(\alpha)^\ell\enspace |\enspace |{\cal I}_{\bar k,r}|>0]\mathbb P(|{\cal I}_{\bar k,r}|>0)   
       + \E[E_{\bar k}(\alpha)^\ell\enspace |\enspace |{\cal I}_{\bar k,r}|=0]\mathbb P(|{\cal I}_{\bar k,r}|=0)\\
       \leq & e^{-\alpha\ell q} + \delta  \frac{e^{-\alpha \ell q}}{N^\ell}\E \lrq{e^{\alpha\ell\max_{i=1:N} f(\xbf^i_k)} \enspace |\enspace |{\cal I}_{\bar k,r}|=0}\\
       \leq& e^{-\alpha\ell q}\lr{1+\delta\frac{1}{N^\ell}e^{\alpha\ell M_f}}
    \end{aligned}
$$
where the last inequality follows from Corollary \ref{boundedfgf}.
Using this inequality in \eqref{newsplit2} we get
\be \begin{aligned}
    \E\lrq{\|\hat x^\alpha_{\bar k}-x^*\|^2}&\leq \frac{1}{4}C_{\bar k} 
    + C_{\bar k}\alpha M_v \ca^{\bar k}D_0^{1/2} + \alpha^2M_v^2\ca^{2\bar k}D_0 \\
    &+\E\lrq{V_{\bar k}^2}^{1/2}e^{-2\alpha q}\lr{N^2+\delta^{1/2}e^{2\alpha M_f}}\\
    &+ C_{\bar k}  \E\lrq{V_{\bar k}}^{1/2}e^{-\alpha q}\lr{N+\delta^{1/2}e^{\alpha M_f}} \\
    & + 2\alpha  M_v \ca^{\bar k}D_0^{1/2}\E\lrq{V_{\bar k}^2}^{1/4} e^{-\alpha q}\lr{N+\delta^{1/4}e^{\alpha M_f}},
\end{aligned}
\ee
and the assumption $\delta \leq\hat\delta e^{-4\alpha M_f}$ yields
\be \begin{aligned}
    \E\lrq{\|\hat x^\alpha_{\bar k}-x^*\|^2}&\leq \frac{1}{4}C_{\bar k} + 
    + C_{\bar k}\alpha M_v \ca^{\bar k}D_0^{1/2} + \alpha^2M_v^2\ca^{2\bar k}D_0 \\
    &+\E\lrq{V_{\bar k}^2}^{1/2}e^{-2\alpha q}\lr{N^2+\hat\delta^{1/2}}\\
    &+ C_{\bar k}  \E\lrq{V_{\bar k}}^{1/2}e^{-\alpha q}\lr{N+\hat\delta^{1/2}e^{-\alpha M_f}} \\
    & + 2\alpha  M_v \ca^{\bar k}D_0^{1/2}\E\lrq{V_{\bar k}^2}^{1/4} e^{-\alpha q}\lr{N+\hat\delta^{1/4}}.
\end{aligned}
\ee
    Proceeding as in the proof of Theorem \ref{thmVkdet}, 
it is easy to see that if $\alpha=\hat \alpha$ sufficiently large  and $M_v$ sufficiently small,  
the inequality above implies that
$\E[\|\hat\xbf^\alpha_{\bar k} - \xbf^*\|^2]\leq C_{\bar k}^2.$ 
\qed\end{proof}

\bibliographystyle{spmpsci}\bibliography{references}

\end{document}